\documentclass[12pt]{article}

\usepackage[all]{xy}

\usepackage[latin1]{inputenc}
\usepackage{amsfonts}
\usepackage{amsmath}
\usepackage{amssymb}
\usepackage{amscd}
\usepackage{amsthm}
\usepackage{indentfirst}
\usepackage[hmargin=2.98cm
,vmargin=4.0cm
]{geometry}

\usepackage{faktor}

\usepackage{tikz-cd}

\usepackage{epsfig}
\usepackage[T1]{fontenc}
\usepackage{hyperref}

\newtheorem{thm}{Theorem}[section]
\newtheorem{lemma}[thm]{Lemma}

\newtheorem{prop}[thm]{Proposition}
\newtheorem{coroll}[thm]{Corollary}
\newtheorem{claim}[thm]{Claim}

\theoremstyle{definition}

\newtheorem{defin}[thm]{Definition}
\newtheorem{rem}[thm]{Remark}
\newtheorem{exam}[thm]{Example}

\newtheorem{question}[thm]{Question}

\newcommand{\R}{{\mathbb{R}}}

\newcommand{\Q}{{\mathbb{Q}}}
\newcommand{\T}{{\mathbb{T}}}
\newcommand{\Z}{{\mathbb{Z}}}
\newcommand{\N}{{\mathbb{N}}}

\newcommand{\CP}{{\mathbb{C}P}}

\newcommand{\cP}{{\mathcal{P}}}

\DeclareMathOperator{\spn}{span}

\DeclareMathOperator{\dist}{dist}

\DeclareMathOperator{\id}{id}
\DeclareMathOperator{\area}{area}

\DeclareMathOperator{\pr}{pr}

\DeclareMathOperator{\supp}{supp}

\DeclareMathOperator{\Ham}{Ham}

\DeclareMathOperator{\Int}{Int}

\DeclareMathOperator{\diam}{diam}

\DeclareMathOperator{\0}{\operatorname{\mathbf{0}}}

\begin{document}
	
	\title{Quantitative thick-fiber property}
	\author{Adi Dickstein\footnote{Partially supported by the Fondation Courtois and the Jenny Panitch Beckow Memorial Scholarship Canada, 2026.}}
	
	\date{}
	
	\setcounter{tocdepth}{2}
	
	\renewcommand{\labelenumi}{(\roman{enumi})}

	\renewcommand{\labelenumi}{(\roman{enumi})}

	\maketitle

\begin{abstract}
We introduce the quantitative thick-fiber property, an analogue of the Entov--Polterovich non-displaceable fiber theorem for maps whose coordinate functions need not Poisson commute. We prove that this property holds for closed symplectic manifolds admitting a suitable, not necessarily linear, functional on the space of continuous functions. Conversely, we establish an obstruction and show that the property fails for certain irrational symplectic four-tori and Kodaira--Thurston manifolds. These results reveal a new rigidity--flexibility distinction: whereas the Entov--Polterovich theorem gives a universal fiber-rigidity phenomenon, its quantitative counterpart holds on some closed symplectic manifolds and fails on others, thereby defining a new symplectic rigidity invariant. We also introduce the non-displaceability width of a continuous map, defined as the infimum of the radii for which a non-displaceable thick fiber exists, and use it to obtain quantitative obstructions to uniform approximation by maps with small Poisson brackets.
\end{abstract}
	\tableofcontents
    
\section{Introduction}\label{s: intro}

Big fiber theorems form a recurring theme in geometry and topology. They assert that a map from a sufficiently rich space to a simpler target must have a large fiber. The meaning of largeness depends on the setting. It may be expressed in terms of dimension, volume, topology, symmetry, or intersections with a distinguished family of subsets.

A classical example is Hurewicz's dimension-lowering theorem, which guarantees that a map from a high-dimensional compact metric space to a lower-dimensional space has a fiber carrying the missing dimension; see \cite[Theorem~VI.7, pp.~91--93]{Hurewicz_Wallman_1948_Dimension_theory}. The Borsuk--Ulam theorem provides a different form of largeness: every continuous map from $S^n$ to $\R^n$ has a fiber containing two antipodal points; see \cite{Borsuk_1933_Drei_Satze}. Thus, in this case, the fiber is large with respect to the symmetry of the domain.

Metric and cohomological versions of this principle were established by Gromov. His waist theorem guarantees a fiber whose neighborhoods have large volume \cite{Gromov_2003_waists}, while his maximal-fiber theorem guarantees a fiber carrying a quantitatively large part of the cohomology of the domain \cite{Gromov_2009_singularities_I,Gromov_2010_singularities_II}. In the topological centerpoint theorem of Karasev, a large fiber is required to intersect every face of a prescribed dimension in a simplex; see \cite{Karasev_2014_Covering_dimension_using_toric_varieties}. Several of these big fiber theorems are discussed from a unified perspective in \cite{DGPZ_2024_Symp_top_and_IVMs}.

The symplectic big-fiber theorem that motivates the present paper is due to Entov--Polterovich. In this theorem, the size of a fiber is measured by its non-displaceability under Hamiltonian diffeomorphisms (see Definition~\ref{def:displaceable}). Entov--Polterovich proved that every smooth map from a closed symplectic manifold to a Euclidean space has a non-displaceable fiber, provided that its coordinate functions pairwise Poisson commute (see Definition~\ref{def:poisson-bracket}).

The assumption in Theorem~\ref{thm:EP-fiber} that the components of $F$ pairwise Poisson commute has a natural interpretation in Hamiltonian mechanics. Poisson-commuting functions arise naturally in Hamiltonian dynamics and play a fundamental role in integrable systems. From this point of view, Theorem~\ref{thm:EP-fiber} asserts that every family of Poisson-commuting functions has a symplectically large joint level set, namely, a non-displaceable one.

Exact Poisson commutativity is, however, an idealized assumption. Suppose that $F=(F_1,\ldots,F_k)\colon M\to\R^k$ has pairwise Poisson-commuting components and that $G=(G_1,\ldots,G_k)\colon M\to\R^k$ is a small $C^1$-perturbation of $F$. Then the Poisson brackets of the components of $G$ are small, but they need not vanish. Consequently, Theorem~\ref{thm:EP-fiber} no longer applies to $G$. In fact, the mere smallness of these Poisson brackets does not guarantee the existence of a non-displaceable fiber: maps with arbitrarily small Poisson brackets may have all their fibers displaceable.

The purpose of this paper is to formulate and study a quantitative stability version of the non-displaceable fiber theorem. To this end, we introduce the quantitative thick-fiber property in Definition~\ref{def: quantitative thick-fiber property}. Instead of asking whether an almost Poisson-commuting map has a non-displaceable fiber, we ask whether it has a non-displaceable \emph{thick fiber}, and how much the fiber must be thickened in terms of the Poisson brackets of the components and the dimension of the target. The desired radius should converge to zero as the Poisson brackets converge uniformly to zero, so that the Entov--Polterovich theorem is recovered in the Poisson-commuting case.

This question has both a positive and a negative aspect. On manifolds carrying a suitable \textit{quasi-state-type functional} (see Definition~\ref{def: ssqs}), we obtain an explicit bound, in terms of the target dimension and the size of the Poisson brackets, on the radius needed to guarantee the existence of a pseudoheavy, and hence non-displaceable, thick fiber; see Theorem~\ref{thm:thick-fiber}. On the other hand, we exhibit closed symplectic manifolds that do not satisfy the quantitative thick-fiber property; see Theorem~\ref{thm: obstruction} and Section~\ref{ss: obst_to_BTFP}.

From the viewpoint of rigidity and flexibility, the Entov--Polterovich theorem gives a universal rigidity statement: on every closed symplectic manifold, exact Poisson commutativity forces the existence of a non-displaceable fiber. The quantitative setting reveals a different picture. On some closed symplectic manifolds, a quantitative remnant of this rigidity persists when exact Poisson commutativity is weakened and exact fibers are replaced by thick fibers. On other closed symplectic manifolds, no estimate of the required kind exists. This distinction is global: the quantitative thick-fiber property does not concern the stability of one fixed map, but asks for an estimate that applies to every smooth map from the given symplectic manifold to every Euclidean space. Thus, unlike the Entov--Polterovich theorem itself, the quantitative thick-fiber property distinguishes between closed symplectic manifolds and serves as a new symplectic rigidity invariant.

\subsection{Definitions and main results}

A \textbf{symplectic manifold} is a pair $(M,\omega)$, where $M$ is a smooth manifold of dimension $2n$ and $\omega$ is a closed $2$-form such that $\omega^n$ is nowhere vanishing. Let $H\colon M\to\R$ be a smooth function. Its \textbf{Hamiltonian vector field} $X_H$ is defined by
\[\omega(X_H,\cdot)=-dH.\]
The flow of $X_H$ is called the \textbf{Hamiltonian flow} of $H$ and is denoted by $(\varphi_H^t)_{t\in I}$, where $I\subseteq\R$ is an interval. Throughout this paper, the manifold $M$ is closed. Therefore, every Hamiltonian vector field on $M$ is complete, and its flow is defined for every $t\in\R$.

More generally, let $H\colon[0,1]\times M\to\R$ be a smooth time-dependent Hamiltonian, and put $H_t=H(t,\cdot)$ for every $t\in[0,1]$. Let $(\varphi_H^t)_{t\in[0,1]}$ denote the Hamiltonian isotopy generated by the time-dependent vector field $X_{H_t}$. A diffeomorphism of $M$ is called a \textbf{Hamiltonian diffeomorphism} if it is the time-$1$ map of a Hamiltonian isotopy. The collection of Hamiltonian diffeomorphisms forms a group, denoted by
\[\Ham(M,\omega)=\{\varphi_H^1:H\in C^\infty([0,1]\times M,\R)\}.\]

\begin{defin}\label{def:poisson-bracket}
The \textbf{Poisson bracket} is the bilinear map
\[\{\cdot,\cdot\}\colon C^\infty(M)\times C^\infty(M)\to C^\infty(M)\]
defined by
\[\{F,G\}=dF(X_G)\]
for every $F,G\in C^\infty(M)$.
\end{defin}

The Poisson bracket is skew-symmetric. Moreover, let $H\colon M\to\R$ be a smooth function, and let $(\varphi_H^t)_{t\in\R}$ be its Hamiltonian flow. Then, for every $F\in C^\infty(M)$, every $x\in M$, and every $t\in\R$, we have
\[\frac{d}{dt}F(\varphi_H^t(x))=\{F,H\}(\varphi_H^t(x)).\]

We also recall the notion of a displaceable set.

\begin{defin}\label{def:displaceable}
A subset $X\subseteq M$ is called \textbf{displaceable} if there exists $\varphi\in\Ham(M,\omega)$ such that
\[\varphi(X)\cap\overline{X}=\varnothing.\]
Otherwise, $X$ is called \textbf{non-displaceable}.
\end{defin}
\begin{exam}
On the sphere $S^2$, a smoothly embedded closed disk $D$ is displaceable if and only if $\area(D)<\frac{1}{2}\area(S^2)$. More generally, a finite disjoint union of smoothly embedded closed disks is displaceable if and only if its total area is strictly smaller than $\frac{1}{2}\area(S^2)$. By contrast, an equator is non-displaceable, although it has zero area.
\end{exam}

For every continuous function $F\colon M\to\R$, we denote by $\supp F$ the closure of its nonzero set. If $M$ is compact, we put $\|F\|_\infty=\max_{x\in M}|F(x)|$.

We now recall the Entov--Polterovich theorem on non-displaceable fibers.

\begin{thm}[Entov--Polterovich, \cite{Entov_Polterovich_2006_QS_and_symp_intersection}]\label{thm:EP-fiber}
Let $(M,\omega)$ be a closed symplectic manifold, and let
\[
F=(F_1,\ldots,F_N)\colon M\longrightarrow\R^N
\]
be a smooth map. Suppose that $\{F_i,F_j\}=0$ for every $i,j\in\{1,\ldots,N\}$. Then there exists $y_0\in F(M)$ such that the fiber $F^{-1}(y_0)$ is non-displaceable.
\end{thm}

It is natural to ask what remains true when the components of $F$ do not Poisson commute exactly, but their Poisson brackets are small. For every integer $k\geq1$ and every smooth map
\[
F=(F_1,\ldots,F_k)\colon M\to\R^k,
\]
define
\[
pb(F)=\max_{1\leq i,j\leq k}\|\{F_i,F_j\}\|_\infty.
\]
Thus, $pb(F)=0$ if and only if the components of $F$ pairwise Poisson commute.

The conclusion of Theorem~\ref{thm:EP-fiber} does not remain true under the weaker assumption that $pb(F)$ is small. Namely, even when all the Poisson brackets are arbitrarily small, it may happen that every fiber of $F$ is displaceable. This leads us to replace an exact fiber by a thick fiber.

For every $y\in\R^k$ and every $r>0$, denote by $B_y(r)$ and $\overline{B}_y(r)$ the Euclidean open and closed balls of radius $r$ centered at $y$, respectively. We also put $\overline{B}_y(0)=\{y\}$. Let $F\colon M\to\R^k$ be a continuous map. For every $y\in\R^k$ and every $r\geq0$, the set
\[
F^{-1}\left(\overline{B}_y(r)\right)
\]
is called the \textbf{closed $r$-thick fiber of $F$ centered at $y$}.

The following example shows why exact fibers must be replaced by thick fibers and illustrates the square-root dependence on the Poisson brackets.

\begin{exam}\label{exam:sphere-thick-fiber}
Equip the unit sphere $S^2\subseteq\R^3$ with the standard area form, whose total area is $4\pi$. Fix $\epsilon>0$, and define the smooth map $F_\epsilon\colon S^2\to\R^2$ by
\[
F_\epsilon(x,y,z)=\left(\epsilon x,2\epsilon y\right)
\]
for every $(x,y,z)\in S^2$. The Poisson bracket of the coordinate functions satisfies $\{x,y\}=-z$. Therefore,
\[
pb(F_\epsilon)
=\left\|\{\epsilon x,2\epsilon y\}\right\|_\infty
=2\epsilon^2.
\]
Since $pb(F_\epsilon)=2\epsilon^2>0$, the hypotheses of Theorem~\ref{thm:EP-fiber} are not satisfied. Moreover, every fiber of $F_\epsilon$ consists of at most two points and is therefore displaceable. Thus, we can arrange that $pb(F_\epsilon)$ is as small as we want, while every fiber of $F_\epsilon$ remains displaceable.

For every $r\geq0$, denote by $X_{\epsilon,r}$ the closed $r$-thick fiber of $F_\epsilon$ centered at the origin $0\in\R^2$, that is,
\[
X_{\epsilon,r}
=F_\epsilon^{-1}\left(\overline{B}_0(r)\right)
=\left\{(x,y,z)\in S^2:\epsilon^2x^2+4\epsilon^2y^2\leq r^2\right\}.
\]
Thus, $X_{\epsilon,r}$ is the intersection of $S^2$ with a solid elliptic cylinder whose axis is the $z$-axis.

Let $E=\{y=0\}\subseteq S^2$ be the equator. We have $E\subseteq X_{\epsilon,r}$ if and only if $r\geq\epsilon$. Therefore, $X_{\epsilon,r}$ is non-displaceable for every $r\geq\epsilon$.

On the other hand, for every $0\leq r<\epsilon$, the set $X_{\epsilon,r}$ is displaceable. Indeed, if $0<r<\epsilon$, then $X_{\epsilon,r}$ is a disjoint union of two smoothly embedded closed disks. Moreover,
\[
X_{\epsilon,r}\subseteq
\left\{|y|\leq\frac{r}{2\epsilon}\right\}.
\]
The spherical band on the right-hand side has area $2\pi r/\epsilon$, and hence
\[
\area(X_{\epsilon,r})
\leq\frac{2\pi r}{\epsilon}
<2\pi.
\]
Therefore, $X_{\epsilon,r}$ is displaceable. By inclusion, the same conclusion holds for $r=0$.

We conclude that $X_{\epsilon,r}$ is non-displaceable if and only if
\[
r\geq\epsilon=\frac{1}{\sqrt{2}}\sqrt{pb(F_\epsilon)}.
\]
Thus, although $F_\epsilon$ has no non-displaceable fiber, its thick fiber centered at the origin becomes non-displaceable at a radius proportional to $\sqrt{pb(F_\epsilon)}$.
\end{exam}

Motivated by this example, we ask whether every almost Poisson-commuting map must have a non-displaceable thick fiber whose radius can be estimated in terms of the target dimension and the size of its Poisson brackets. The following definition formalizes this question.

\begin{defin}\label{def: quantitative thick-fiber property}
Let $(M,\omega)$ be a closed symplectic manifold. Let
\[
f\colon\N\to(0,\infty),\qquad
g\colon[0,\infty)\to[0,\infty)
\]
be functions such that $g$ is continuous and $g(0)=0$.

We say that $(M,\omega)$ has the \textbf{quantitative thick-fiber property with respect to $(f,g)$} if, for every integer $k\geq1$ and every smooth map $F\colon M\to\R^k$, there exists $y\in F(M)$ such that the closed $f(k)\cdot g\left(pb(F)\right)$-thick fiber of $F$ centered at $y$, namely,
\[
F^{-1}\left(
\overline{B}_y\left(f(k)g\left(pb(F)\right)\right)
\right),
\]
is non-displaceable.

We say that $(M,\omega)$ has the \textbf{quantitative thick-fiber property} if there exist functions $f$ and $g$ as above such that $(M,\omega)$ has the quantitative thick-fiber property with respect to $(f,g)$.
\end{defin}
The function $f$ records the dependence of the radius on the dimension of the target, while $g$ records its dependence on the failure of the components of $F$ to Poisson commute. The condition $g(0)=0$ ensures that the radius vanishes when the components of $F$ pairwise Poisson commute.

The original proof of Theorem~\ref{thm:EP-fiber} uses the notion of a partial symplectic quasi-state, which is a not necessarily linear functional $\zeta\colon C(M)\to\R$ satisfying several axioms. Here we define a similar concept that will be useful in the formulation of our main result.

\begin{defin}\label{def: ssqs}
Let $(M,\omega)$ be a closed symplectic manifold and let $C>0$. A functional $\zeta\colon C(M)\to\R$ is called a \textbf{semi-symplectic quasi-state with constant $C$} if it satisfies the following axioms:
\begin{itemize}
    \item (Normalization) $\zeta(c)=c$ for every $c\in\R$.
    \item (Monotonicity) For every $F,G\in C(M)$, if $F\leq G$, then $\zeta(F)\leq\zeta(G)$.
    \item (Vanishing) For every $F\in C(M)$, if $F$ has displaceable support, then $\zeta(F)=0$.
    \item (One-sided Poisson-bracket inequality) For every $F,G\in C^\infty(M)$,
    \[\zeta(F+G)-\zeta(F)-\zeta(G)\leq C\sqrt{\|\{F,G\}\|_\infty}.\]
\end{itemize}
\end{defin}

In Section~\ref{ss:semi-symp_qs-examples}, we compare this notion with existing quasi-state-type functionals.

We also adapt to semi-symplectic quasi-states the notion of a pseudoheavy set, introduced by Kawasaki--Orita for partial symplectic quasi-states in \cite{Kawasaki_Orita_2021_seudoheavy}. See Section~\ref{ss: superheavy thick-fibers} for further discussion.

\begin{defin}\label{def: pseudoheavy}
Let $(M,\omega)$ be a closed symplectic manifold and let $\zeta$ be a semi-symplectic quasi-state on $(M,\omega)$. A closed subset $X\subseteq M$ is called \textbf{$\zeta$-pseudoheavy} if, for every open neighborhood $U$ of $X$, there exists a continuous function $F\in C(M)$ supported in $U$ such that $\zeta(F)>0$.
\end{defin}

\begin{rem}\label{rem: pseudoheavy implies non-disp}
Kawasaki--Orita proved that every pseudoheavy set with respect to a partial symplectic quasi-state is non-displaceable; see \cite[Proposition~1.6]{Kawasaki_Orita_2021_seudoheavy}. The same conclusion holds for semi-symplectic quasi-states. Indeed, if a closed subset $X\subseteq M$ were displaceable, then it would admit a displaceable open neighborhood $U$. The vanishing property would imply that $\zeta(F)=0$ for every continuous function $F$ supported in $U$, contradicting the pseudoheaviness of $X$.
\end{rem}

We are now ready to state our main thick-fiber theorem.

\begin{thm}\label{thm:thick-fiber}
Let $(M,\omega)$ be a closed symplectic manifold, let $C>0$, and assume that $(M,\omega)$ is equipped with a semi-symplectic quasi-state $\zeta$ with constant $C$. Then, for every integer $k\geq1$ and every smooth map $F\colon M\to\R^k$, there exists $y_0\in F(M)\subset \R^k$ such that
\[F^{-1}\left(\overline{B}_{y_0}\left(4C\cdot k^{3/2}\cdot\sqrt{pb(F)}\right)\right)\]
is $\zeta$-pseudoheavy, and in particular, non-displaceable.
\end{thm}

Examples of closed symplectic manifolds admitting semi-symplectic quasi-states include closed oriented surfaces, complex projective spaces, and products of spheres; see Section~\ref{ss:semi-symp_qs-examples} for details. Consequently, all these manifolds satisfy the quantitative thick-fiber property by the following immediate corollary of Theorem~\ref{thm:thick-fiber}.

\begin{coroll}\label{cor:main-2}
Let $(M,\omega)$ be a closed symplectic manifold, let $C>0$, and assume that $(M,\omega)$ is equipped with a semi-symplectic quasi-state with constant $C$. Define the functions $f\colon\N\to(0,\infty)$ and $g\colon[0,\infty)\to[0,\infty)$ by
\[f(k)=4Ck^{3/2},\qquad g(t)=\sqrt{t}\]
for every $k\in\N$ and every $t\in[0,\infty)$. Then $(M,\omega)$ has the quantitative thick-fiber property with respect to the functions $f$ and $g$.
\end{coroll}
\begin{rem}
The proof of Theorem~\ref{thm:thick-fiber} uses a special partition of unity on $M$, whose elements are divided into $k+1$ colors, where $k$ is the dimension of the target of the map. The functions assigned the same color have pairwise disjoint positivity sets, so each color class can be combined into a single function. Thus, the main quantitative argument involves only $k+1$ functions. Suitable estimates on the Poisson brackets of these color sums, together with the one-sided Poisson-bracket inequality, imply that one of the original partition functions has a positive $\zeta$-value. Its support is contained in a thick fiber of controlled radius, and a limiting argument produces a pseudoheavy thick fiber. The vanishing property ensures that every pseudoheavy set is non-displaceable. The use of $k+1$ colors was inspired by the proof of \cite[Theorem~6.1]{DGPZ_2024_Symp_top_and_IVMs}. The complete proof is given in Section~\ref{ss: proof of thick-fiber}.
\end{rem}
\begin{rem}
After the main quantitative thick-fiber result of this paper had been obtained, the author learned of related work by Jun Zhang in contact geometry, called ``Contact rigidity and Poisson bracket invariants.'' The author and Zhang then exchanged preliminary drafts. At the time of this exchange, both the symplectic result of the present paper and Zhang's quantitative contact big-fiber theorem, which generalizes the contact big-fiber theorem proved in \cite{SUV_2026_contact_big_fiber_thm,UZ_2025_contact_qs,DUZ_2026_quantitative_cntct_ham_dynmcs}, had already been obtained independently.

The two arguments share the general strategy of decomposing the target into sets of controlled diameter, pulling the resulting functions back by the given map, and using a quantitative Poisson-bracket estimate together with a vanishing property to obtain a non-displaceable thick fiber. The settings and the quantitative constructions are different. Zhang's theorem concerns contact displaceability and involves both contact brackets and Reeb derivatives. The proof presented here uses the Freudenthal--Kuhn triangulation and a $(k+1)$-coloring and produces a pseudoheavy thick fiber with radius of order $k^{3/2}\sqrt{pb(F)}$.
\end{rem}

It is interesting to understand whether the functions $f$ and $g$ from Corollary~\ref{cor:main-2} are optimal. The following result provides lower bounds on the growth of both functions and, in particular, shows that the choice of $g$ in Corollary~\ref{cor:main-2} is optimal up to multiplication by a positive constant.

\begin{thm}\label{thm:lower-bound-thick-fiber-radius}
Let $(M,\omega)$ be a closed symplectic manifold of positive dimension. Let $f\colon\N\to(0,\infty)$ and $g\colon[0,\infty)\to[0,\infty)$ be functions, where $g$ is continuous and satisfies $g(0)=0$. Suppose that $(M,\omega)$ has the quantitative thick-fiber property with respect to $f$ and $g$. Then there exist constants $c_f,c_g>0$ such that, for every integer $k\geq1$ and every $t\geq0$, we have
\[
f(k)\geq c_f\sqrt{k},\qquad\text{and}\qquad g(t)\geq c_g\sqrt{t}.
\]
\end{thm}
The proof of Theorem~\ref{thm:lower-bound-thick-fiber-radius} is given in Section~\ref{ss: Proof of the lower-bound-thick-fiber-radius}.
\begin{rem}
The question of whether the growth rate $k^{3/2}$ of the function $f$ from Corollary~\ref{cor:main-2} is optimal remains open. Let us mention that, in the proof of Theorem~\ref{thm:thick-fiber}, we obtain this function using a refined partition-of-unity construction. More naive constructions, which we used at the beginning of this project, yield a function $f$ that grows exponentially with the dimension.
\end{rem}

We next turn to obstructions to the quantitative thick-fiber property. By negating Definition~\ref{def: quantitative thick-fiber property}, a closed symplectic manifold $(M,\omega)$ does not have this property if, for every pair of functions
\[
f\colon\N\to(0,\infty),\qquad
g\colon[0,\infty)\to[0,\infty),
\]
where $g$ is continuous and satisfies $g(0)=0$, there exist an integer $k\geq1$ and a smooth map $F\colon M\to\R^k$ such that, for every $y\in F(M)$, the closed $f(k)g\left(pb(F)\right)$-thick fiber of $F$ centered at $y$, namely,
\[
F^{-1}\left(
\overline{B}_y\left(f(k)g\left(pb(F)\right)\right)
\right),
\]
is displaceable.

Our next aim is to give a sufficient condition for a closed symplectic manifold not to have the quantitative thick-fiber property. To state the result, we introduce the Poisson bracket of circle-valued functions. Identify $S^1=\R/\Z$, denote by $q\colon\R\to S^1$ the quotient map, and let $\theta\in\Omega^1(S^1)$ be the unique smooth one-form satisfying $q^*\theta=dt$. Let $\Pi$ be the Poisson bivector associated with $\omega$. For every pair of smooth maps $u,v\colon M\to S^1$, define
\[
\{u,v\}=\Pi(u^*\theta,v^*\theta).
\]
Equivalently, if $\widetilde u,\widetilde v\colon U\to\R$ are local lifts of $u$ and $v$, then
\[
\{u,v\}|_U=\{\widetilde u,\widetilde v\}.
\]
In particular, the definition does not depend on the choice of the local lifts.

\begin{thm}\label{thm: obstruction}
Let $(M,\omega)$ be a closed symplectic manifold. Suppose that there exist smooth maps $a,b,c\colon M\to S^1$ and constants $A,B\in\R$ such that
\[\{a,c\}=A,\qquad \{b,c\}=B.\]
If $B\neq0$ and $A/B\notin\Q$, then $(M,\omega)$ does not have the quantitative thick-fiber property.
\end{thm}

The proof of Theorem~\ref{thm: obstruction} is given in Section~\ref{ss: Proof of obstruction}. In Section~\ref{ss: obst_to_BTFP}, we present several examples of closed symplectic manifolds satisfying the hypotheses of Theorem~\ref{thm: obstruction}, and therefore not having the quantitative thick-fiber property. These examples include irrational symplectic four-tori and the Kodaira--Thurston manifold equipped with certain symplectic forms.

Combining Corollary~\ref{cor:main-2} with Theorem~\ref{thm: obstruction} gives the following obstruction to the existence of semi-symplectic quasi-states.
\begin{coroll}
Let $(M,\omega)$ be a closed symplectic manifold. Suppose that there exist smooth maps $a,b,c\colon M\to S^1$ and constants $A,B\in\R$ such that
\[\{a,c\}=A,\qquad \{b,c\}=B.\]
If $B\neq0$ and $A/B\notin\Q$, then $(M,\omega)$ does not admit a semi-symplectic quasi-state.
\end{coroll}

This corollary is related to the question posed by Entov--Polterovich in \cite{EP_2008_symp_qs_and_semi_simplicity} of whether the standard symplectic $4$-torus, or more generally any closed symplectic manifold of dimension greater than $2$ satisfying $\left.[\omega]\right|_{\pi_2(M)}=0$, admits a nonlinear genuine symplectic quasi-state. As far as we know, this question remains open. Our corollary does not settle this question, because the Poisson-bracket inequality is an additional property that is not included in the axioms of a genuine symplectic quasi-state. However, it obstructs certain closed symplectic manifolds, including the irrational symplectic $4$-tori considered below, from admitting a genuine symplectic quasi-state that satisfies the Poisson-bracket inequality.

\paragraph*{Organization of the paper.}
The remainder of the Introduction is organized as follows. In Section~\ref{ss:semi-symp_qs-examples}, we present examples of semi-symplectic quasi-states, including partial symplectic quasi-states on closed oriented surfaces. In Section~\ref{ss: superheavy thick-fibers}, we discuss conditions under which the pseudoheavy thick fiber provided by Theorem~\ref{thm:thick-fiber} is superheavy. Next, in Section~\ref{ss: non-disp width}, we introduce the non-displaceability width and use it to obtain quantitative lower bounds for the Poisson brackets of smooth maps approximating a given continuous map. In Section~\ref{ss: obst_to_BTFP}, we present examples of symplectic manifolds that do not have the quantitative thick-fiber property. The Introduction concludes in Section~\ref{ss: discussion} with a discussion of Riemannian targets and the displacement energy of fibers. Section~\ref{s: proofs} contains the proofs. We first prove the result concerning partial symplectic quasi-states on surfaces and the lower bounds for the functions appearing in the quantitative thick-fiber property. We then prove the obstruction theorem and, finally, Theorem~\ref{thm:thick-fiber}; its proof includes the construction of the special partition of unity used in the argument.

\subsection{Semi-symplectic quasi-states}\label{ss:semi-symp_qs-examples}

The concept of quasi-states arose from discussions surrounding the von Neumann axioms of quantum mechanics; see \cite[Section 10]{Entov_Polterovich_2006_QS_and_symp_intersection}, \cite[Section 1]{Butler_2022_semisolid_sets}, and the references therein. In 1991, Aarnes adapted this concept to topology and functional analysis, introduced \textit{topological quasi-states}, which are not necessarily linear functionals on the space of real-valued continuous functions on a compact Hausdorff space $X$, and formulated a geometric classification theorem; see \cite{Aarnes_1991_quasi_states_and_quasi_measures}. Topological quasi-states are linear on singly generated subalgebras of $C(X)$; this property is called \textbf{quasi-linearity}.

In 2006, Entov--Polterovich used Floer theory to construct certain not necessarily linear functionals on the space of real-valued continuous functions on closed symplectic manifolds, called \textit{partial symplectic quasi-states}. Under certain assumptions on the symplectic manifold, these partial quasi-states are also topological quasi-states; in such cases, we call them \textit{genuine symplectic quasi-states}. These genuine symplectic quasi-states satisfy a strong version of the quasi-linearity property: each is linear on every Poisson-commuting subalgebra of $C^\infty(M)\subseteq C(M)$.

Given a closed symplectic manifold $(M,\omega)$, a partial symplectic quasi-state $\zeta\colon C(M)\to\R$ satisfies several axioms, including the following:

\begin{itemize}
    \item (Normalization) $\zeta(c)=c$ for every $c\in\R$.
    \item (Monotonicity) For every $F,G\in C(M)$, if $F\leq G$, then $\zeta(F)\leq\zeta(G)$.
    \item (Vanishing) For every $F\in C(M)$, if $F$ has displaceable support, then $\zeta(F)=0$.
    \item (Quasi-subadditivity) For every $F,G\in C^\infty(M)$, if $\{F,G\}=0$, then
    \[\zeta(F+G)\leq\zeta(F)+\zeta(G).\]
\end{itemize}

As mentioned above, Entov--Polterovich proved in \cite{Entov_Polterovich_2006_QS_and_symp_intersection} that, under certain assumptions, their construction of a partial symplectic quasi-state produces a genuine symplectic quasi-state, which also satisfies the following axiom:
\begin{itemize}
    \item (Quasi-linearity) For every $F,G\in C^\infty(M)$, if $\{F,G\}=0$, then the restriction of $\zeta$ to $\spn(F,G)$ is linear.
\end{itemize}
In \cite{Entov_Polterovich_Zapolsky_2007_quasi-morphisms_and_PB}, Entov--Polterovich--Zapolsky proved, under the same assumptions, a quantitative version of the quasi-linearity property in the form of the \textbf{Poisson-bracket inequality}: there exists a constant $C>0$, depending on $\zeta$, such that, for every $F,G\in C^\infty(M)$,
\begin{equation}\label{eq: pb_ineq}
    \left|\zeta(F+G)-\zeta(F)-\zeta(G)\right|\leq C\sqrt{\|\{F,G\}\|_\infty}.
\end{equation}
It is interesting to ask whether partial symplectic quasi-states satisfy a similar inequality. Since quasi-subadditivity implies that
\[\zeta(F+G)-\zeta(F)-\zeta(G)\leq0\]
whenever $\{F,G\}=0$, it is natural to consider a \textbf{one-sided Poisson-bracket inequality}: there exists a constant $C>0$, depending only on the partial symplectic quasi-state $\zeta$, such that, for every $F,G\in C^\infty(M)$,
\[\zeta(F+G)-\zeta(F)-\zeta(G)\leq C\sqrt{\|\{F,G\}\|_\infty}.\]
This one-sided inequality motivates the notion of a semi-symplectic quasi-state introduced in Definition~\ref{def: ssqs}.

In the rest of this section, we discuss several examples of semi-symplectic quasi-states. We start with a two-dimensional example, the \textit{median quasi-state}. It is hard to pinpoint who exactly introduced it, but it is safe to say that it was known to specialists toward the end of the 1990s.

\begin{exam}\label{exam:med-qs}
The median quasi-state is a genuine symplectic quasi-state $\zeta$ on $S^2$ equipped with the normalized standard area form. It can be described as follows. Let $F\colon S^2\to\R$ be a Morse function. Then $\zeta(F)$ is the unique value of $F$ whose level set contains a connected component $C$ such that every connected component of $S^2\setminus C$ has area at most $\frac{1}{2}$. Such a value always exists.

The functional $\zeta$ is $1$-Lipschitz with respect to the uniform norm on the space of Morse functions. Therefore, it has a unique continuous extension to the whole space $C(S^2)$.

It was proved in \cite{Entov_Polterovich_Zapolsky_2007_quasi-morphisms_and_PB} that the median quasi-state satisfies the Poisson-bracket inequality. More recently, it was proved in \cite{HK-K-K-R_2026_J-triangles_and_med_qs} that the optimal constant in the Poisson-bracket inequality for the median quasi-state is $\frac{1}{2}$.
\end{exam}

Let us review the Entov--Polterovich construction, assuming that the reader is familiar with the notion of spectral invariants in Floer theory. See \cite{Entov_Polterovich_2003_Calabi_quasimorphism_and_quantum_homology,Oh_1997_Symplectic_topology_I,Oh_1999_Symplectic_topology_II,Oh_2005_Construction_of_spectral_invariants,Schwarz_action_spectrum_for_aspherical_mfds} for relevant background. Let $(M,\omega)$ be a closed symplectic manifold, and let $e$ be a nonzero idempotent element in the quantum homology $QH_*(M,\omega)$ of $(M,\omega)$, that is, $e*e=e$. Denote by $c(e,\cdot)\colon C^\infty(S^1\times M)\to\R$ the spectral invariant associated with $e$. For every autonomous Hamiltonian $F\in C^\infty(M)$, define
\[\zeta_e(F)=\lim_{k\to\infty}\frac{c(e,kF)}{k}.\]
The functional $\zeta_e\colon C^\infty(M)\to\R$ is $1$-Lipschitz with respect to the uniform norm on $C^\infty(M)$. Therefore, it extends uniquely to a functional on $C(M)$.

Entov--Polterovich proved that $\zeta_e$ is a partial symplectic quasi-state. Moreover, if $e$ is the unit of a field factor of $QH_*(M,\omega)$, then $\zeta_e$ is a genuine symplectic quasi-state. Basic examples of symplectic manifolds whose quantum homology has a field factor include $\CP^n$ and products of spheres. On $S^2$, their construction recovers the median quasi-state from Example~\ref{exam:med-qs}. There are also closed symplectic manifolds, such as $S^2\times S^2$, whose quantum homology contains several field factors. Thus, the Entov--Polterovich construction produces, a priori, several different genuine symplectic quasi-states on such manifolds. See \cite{Eliashberg_Polterovich_2010_qs} for further discussion.

Entov--Polterovich--Zapolsky proved in \cite{Entov_Polterovich_Zapolsky_2007_quasi-morphisms_and_PB} that, if $e$ is the unit of a field factor of $QH_*(M,\omega)$, then $\zeta_e$ satisfies the Poisson-bracket inequality~\eqref{eq: pb_ineq} for some constant $C>0$. It follows immediately from the definition that a genuine symplectic quasi-state satisfying the Poisson-bracket inequality is, in particular, a semi-symplectic quasi-state.

\begin{exam}
Let $(M,\omega)$ be a closed symplectic manifold, let $k\geq1$, and let $\zeta_1,\ldots,\zeta_k$ be semi-symplectic quasi-states on $(M,\omega)$ with constants $C_1,\ldots,C_k>0$. Define the functional $\zeta\colon C(M)\to\R$ by
\[\zeta(F)=\max_{1\leq i\leq k}\zeta_i(F)\]
for every $F\in C(M)$, and put
\[C=\max_{1\leq i\leq k}C_i.\]
Let us show that $\zeta$ is a semi-symplectic quasi-state with constant $C$. The normalization, monotonicity, and vanishing properties of $\zeta$ follow immediately from the definition.

Let $F,G\in C^\infty(M)$. There exists $i\in\{1,\ldots,k\}$ such that $\zeta(F+G)=\zeta_i(F+G)$. By the definition of $\zeta$, we have $\zeta_i(F)\leq\zeta(F)$ and $\zeta_i(G)\leq\zeta(G)$. Therefore,
\[\zeta(F+G)-\zeta(F)-\zeta(G)\leq\zeta_i(F+G)-\zeta_i(F)-\zeta_i(G)\leq C_i\sqrt{\|\{F,G\}\|_\infty}\leq C\sqrt{\|\{F,G\}\|_\infty}.\]
Thus, $\zeta$ satisfies the one-sided Poisson-bracket inequality with constant $C$.
\end{exam}

A nontrivial example of a semi-symplectic quasi-state is provided by a partial symplectic quasi-state on a surface.

\begin{prop}\label{prop:ssqs_on_surfaces}
Let $(\Sigma,\sigma)$ be a closed oriented surface equipped with an area form $\sigma$, and put $A=\int_\Sigma\sigma$. Then every partial symplectic quasi-state on $(\Sigma,\sigma)$ is a semi-symplectic quasi-state with constant $2\sqrt{2A}$.
\end{prop}
The proof of Proposition~\ref{prop:ssqs_on_surfaces} is given in Section~\ref{ss: proof of prop: ssqs_on_surfaces}.
\begin{rem}
If $\Sigma$ has positive genus, then the Entov--Polterovich construction produces a partial symplectic quasi-state on $\Sigma$ that is not genuine; see, for instance, \cite{HLS_2016_Spec_inv_in_2D}. Thus, such quasi-states provide nontrivial examples of semi-symplectic quasi-states that are not genuine.
\end{rem}

\subsection{Superheavy thick fibers}\label{ss: superheavy thick-fibers}

Kawasaki--Orita introduced pseudoheavy sets with respect to partial symplectic quasi-states in \cite{Kawasaki_Orita_2021_seudoheavy}. This notion fits into the hierarchy of heavy and superheavy sets introduced by Entov--Polterovich in \cite{EP_2009_rigid_subsets}. In the same paper, Kawasaki--Orita strengthened the Entov--Polterovich fiber theorem by proving the existence of a pseudoheavy fiber with respect to a partial symplectic quasi-state. Later, Mak--Sun--Varolg\"une\c{s} further strengthened this conclusion for partial symplectic quasi-states arising from the Entov--Polterovich construction by proving the existence of a heavy fiber; see \cite{MSV_2024_heavy_sets_and_SH}.

In our setting, the corresponding superheavy conclusion follows readily from Theorem~\ref{thm:thick-fiber} and existing results. We now recall the notions of heavy, superheavy, and pseudoheavy subsets and the relationships between them that yield this conclusion.

\begin{defin}
Let $(M,\omega)$ be a closed symplectic manifold equipped with a partial symplectic quasi-state $\zeta$, and let $X\subseteq M$ be a nonempty closed subset.
\begin{itemize}
\item (\cite{EP_2009_rigid_subsets}) The set $X$ is called \textbf{$\zeta$-heavy} if, for every continuous function $H\in C(M)$, we have $\zeta(H)\geq\min_X H$.
\item (\cite{EP_2009_rigid_subsets}) The set $X$ is called \textbf{$\zeta$-superheavy} if, for every continuous function $H\in C(M)$, we have $\zeta(H)\leq\max_X H$.
\item (\cite{Kawasaki_Orita_2021_seudoheavy}) The set $X$ is called \textbf{$\zeta$-pseudoheavy} if, for every open neighborhood $U$ of $X$, there exists a continuous function $H\in C(M)$ supported in $U$ and satisfying $\zeta(H)>0$.
\end{itemize}
\end{defin}
Entov--Polterovich proved in \cite{EP_2009_rigid_subsets} that every superheavy set is heavy, that every heavy set is non-displaceable, and that every superheavy set intersects every heavy set. They also proved that, if $\zeta$ is a genuine symplectic quasi-state, then every heavy set is superheavy.
Kawasaki--Orita proved in \cite{Kawasaki_Orita_2021_seudoheavy} that every heavy set is pseudoheavy, that every pseudoheavy set is non-displaceable, and that every superheavy set intersects every pseudoheavy set. Moreover, under the assumption of \textit{simplicity} of the partial symplectic quasi-state, they proved that every pseudoheavy set is heavy.

For a partial symplectic quasi-state $\zeta$, one can define an associated set function $\tau_\zeta$ on the compact subsets of $M$. For every compact subset $K\subseteq M$, define
\[\tau_\zeta(K)=\inf\{\zeta(a):a\in C(M,[0,1]),\ a|_K\equiv1\}.\]
The partial symplectic quasi-state $\zeta$ is called \textbf{simple} if, for every compact subset $K\subseteq M$, we have either $\tau_\zeta(K)=0$ or $\tau_\zeta(K)=1$.\footnote{Aarnes defined in \cite{Aarnes_1993_simple_and_extremal} a different notion of simplicity for quasi-states. For genuine symplectic quasi-states, this is equivalent to the definition presented here; see \cite[Theorem~2.3]{Aarnes_1993_simple_and_extremal} and \cite[Remark~4.3]{Kawasaki_Orita_2021_seudoheavy}.}
Thus, by \cite[Proposition~1.9]{Kawasaki_Orita_2021_seudoheavy}, if $\zeta$ is a simple partial symplectic quasi-state, then every pseudoheavy set is heavy. Combining this with the facts above, we deduce that, if $\zeta$ is a simple genuine symplectic quasi-state, then every pseudoheavy set is superheavy.

Moreover, it was proved in \cite[Corollary~1.10]{MSV_2024_heavy_sets_and_SH} that every genuine symplectic quasi-state obtained from the Entov--Polterovich construction is simple.\footnote{For comparison, Kawasaki--Orita state that the partial symplectic quasi-state on a closed oriented surface obtained from the Entov--Polterovich construction is simple in the sense defined above; see \cite[pp.~3--4]{Kawasaki_Orita_2021_seudoheavy}. In higher dimensions, this question remains open. However, when the surface has positive genus, it was proved in \cite[Section~1.3]{HLS_2016_Spec_inv_in_2D} that the partial symplectic quasi-state obtained from the Entov--Polterovich construction does not satisfy the condition that, for genuine quasi-states, is equivalent to simplicity in the sense of Aarnes.} As mentioned before, every genuine symplectic quasi-state obtained from the Entov--Polterovich construction also satisfies the Poisson-bracket inequality; see \cite{Entov_Polterovich_Zapolsky_2007_quasi-morphisms_and_PB}. This leads to the following consequence of Theorem~\ref{thm:thick-fiber}.

\begin{coroll}
Let $(M,\omega)$ be a closed symplectic manifold, let $C>0$, and assume that $(M,\omega)$ is equipped with a simple genuine symplectic quasi-state $\zeta$ satisfying the Poisson-bracket inequality with constant $C$. Then, for every integer $k\geq1$ and every smooth map $F\colon M\to\R^k$, there exists $y_0\in F(M)$ such that
\[F^{-1}\left(\overline{B}_{y_0}\left(4C\cdot k^{3/2}\cdot\sqrt{pb(F)}\right)\right)\]
is $\zeta$-superheavy.
\end{coroll}
\begin{proof}
Since $\zeta$ satisfies the Poisson-bracket inequality with constant $C$, it is a semi-symplectic quasi-state with constant $C$. By Theorem~\ref{thm:thick-fiber}, there exists $y_0\in F(M)$ such that $F^{-1}\left(\overline{B}_{y_0}\left(4C\cdot k^{3/2}\cdot\sqrt{pb(F)}\right)\right)$ is $\zeta$-pseudoheavy. Since $\zeta$ is simple, this set is $\zeta$-heavy by \cite[Proposition~1.9]{Kawasaki_Orita_2021_seudoheavy}. Finally, since $\zeta$ is genuine, this set is $\zeta$-superheavy, as required.
\end{proof}

\subsection{Non-displaceability width}\label{ss: non-disp width}

Let $(M,\omega)$ be a closed symplectic manifold. For every integer $k\geq1$ and every continuous map $F\colon M\to\R^k$, put
\[
\|F\|_\infty=\max_{x\in M}\|F(x)\|,
\]
where $\|\cdot\|$ denotes the Euclidean norm on $\R^k$.
\begin{defin}
Let $F\colon M\to\R^k$ be a continuous map. We define the \textbf{non-displaceability width} of $F$ by
\[w(F):=\inf\left\{r>0:F^{-1}\left(\overline{B}_p(r)\right)\text{ is non-displaceable for some }p\in\R^k\right\}.\]
\end{defin}
The compactness of $M$ guarantees that the non-displaceability width is finite.

If $(M,\omega)$ has the quantitative thick-fiber property with respect to functions $f\colon\N\to(0,\infty)$ and $g\colon[0,\infty)\to[0,\infty)$, then, by definition, every smooth map $F\colon M\to\R^k$ satisfies
\begin{equation}\label{eq: upper_bound_for_w(F)}
    w(F)\leq f(k)\cdot g\left(pb(F)\right).
\end{equation}
In particular, if $(M,\omega)$ is equipped with a semi-symplectic quasi-state with constant $C>0$, then Theorem~\ref{thm:thick-fiber} implies that every smooth map $F\colon M\to\R^k$ satisfies
\[w(F)\leq4Ck^{3/2}\cdot\sqrt{pb(F)}.\]

The non-displaceability width satisfies the following basic properties.
\begin{prop}\label{prop:width-stability}
Let $F,G\colon M\to\R^k$ be continuous maps. Then
\[\left| w(F)-w(G)\right|\leq\|F-G\|_\infty.\]
Moreover, for every $\lambda>0$,
\[w(\lambda F)=\lambda w(F).\]
\end{prop}
\begin{proof}
Put $\delta=\|F-G\|_\infty$. First, let us show that
\[w(G)\leq w(F)+\delta.\]
Let $p\in\R^k$ and $r>0$, and assume that $F^{-1}\left(\overline{B}_p(r)\right)$ is non-displaceable. If $x\in F^{-1}\left(\overline{B}_p(r)\right)$, then
\[\| G(x)-p\|\leq\| G(x)-F(x)\|+\| F(x)-p\|\leq\delta+r.\]
Consequently,
\[F^{-1}\left(\overline{B}_p(r)\right)\subseteq G^{-1}\left(\overline{B}_p(r+\delta)\right).\]
Since every set containing a non-displaceable set is non-displaceable, the set $G^{-1}\left(\overline{B}_p(r+\delta)\right)$ is non-displaceable. Hence, $w(G)\leq r+\delta$. Taking the infimum over all such $p$ and $r$, we obtain $w(G)\leq w(F)+\delta$. By interchanging $F$ and $G$, we also obtain $w(F)\leq w(G)+\delta$. Consequently, $\left| w(F)-w(G)\right|\leq\delta$.

For the second assertion, note that for every $p\in\R^k$ and every $r>0$,
\[(\lambda F)^{-1}\left(\overline{B}_p(r)\right)=F^{-1}\left(\overline{B}_{p/\lambda}\left(\frac{r}{\lambda}\right)\right).\]
Since $p/\lambda$ ranges over $\R^k$ and $r/\lambda$ ranges over $(0,\infty)$, it follows directly from the definition that $w(\lambda F)=\lambda w(F)$.
\end{proof}

Proposition~\ref{prop:width-stability} will be useful in the proof of Theorem~\ref{thm:lower-bound-thick-fiber-radius} and will also have applications to $C^0$-approximation problems involving the Poisson bracket, as we will see below.

\subsubsection*{$C^0$-approximations and the Poisson bracket}

Assume that $(M,\omega)$ is equipped with a semi-symplectic quasi-state with constant $C>0$. The stability property above gives the following quantitative obstruction to almost Poisson commutativity. For every $t\in\R$, put $t_+:=\max\{t,0\}$.
\begin{coroll}\label{cor:width-pb-obstruction}
Let $F\colon M\to\R^k$ be a continuous map and let $G\colon M\to\R^k$ be a smooth map. Then
\[pb(G)\geq\frac{\left(w(F)-\|F-G\|_\infty\right)_+^2}{16C^2k^3}.\]
\end{coroll}
\begin{proof}
By Proposition~\ref{prop:width-stability}, we have
\[w(F)\leq w(G)+\|F-G\|_\infty.\]
On the other hand, Theorem~\ref{thm:thick-fiber} gives
\[w(G)\leq4Ck^{3/2}\cdot\sqrt{pb(G)}.\]
Therefore,
\[\left(w(F)-\|F-G\|_\infty\right)_+\leq4Ck^{3/2}\cdot\sqrt{pb(G)}.\]
Squaring both sides completes the proof.
\end{proof}

In particular, if $w(F)>0$, then every smooth map $G\colon M\to\R^k$ satisfying $\|F-G\|_\infty<w(F)$ also satisfies
\[pb(G)\geq\frac{\left(w(F)-\|F-G\|_\infty\right)^2}{16C^2k^3}>0.\]

For every continuous map $F\colon M\to\R^k$, define
\[\liminf_{\substack{G\xrightarrow{C^0}F\\G\in C^\infty(M,\R^k)}}pb(G):=\lim_{\epsilon\to0^+}\inf\left\{pb(G):G\in C^\infty(M,\R^k),\ \|F-G\|_\infty<\epsilon\right\}.\]

\begin{coroll}\label{cor:width-liminf}
For every continuous map $F\colon M\to\R^k$, we have
\[\liminf_{\substack{G\xrightarrow{C^0}F\\G\in C^\infty(M,\R^k)}}pb(G)\geq\frac{w(F)^2}{16C^2k^3}.\]
In particular, if $w(F)>0$, then
\[\liminf_{\substack{G\xrightarrow{C^0}F\\G\in C^\infty(M,\R^k)}}pb(G)>0.\]
\end{coroll}
\begin{proof}
Let $\epsilon>0$. By Corollary~\ref{cor:width-pb-obstruction}, every smooth map $G\colon M\to\R^k$ satisfying $\|F-G\|_\infty<\epsilon$ also satisfies
\[pb(G)\geq\frac{\left(w(F)-\|F-G\|_\infty\right)_+^2}{16C^2k^3}\geq\frac{\left(w(F)-\epsilon\right)_+^2}{16C^2k^3}.\]
Therefore,
\[\inf\left\{pb(G):G\in C^\infty(M,\R^k),\ \|F-G\|_\infty<\epsilon\right\}\geq\frac{\left(w(F)-\epsilon\right)_+^2}{16C^2k^3}.\]
Since $\epsilon>0$ was arbitrary, the required inequality follows.
\end{proof}

When $F$ is smooth, $k=2$, and $w(F)>0$, the positivity of the limit inferior also follows from the $C^0$-rigidity of Poisson commutativity proved in \cite[Theorem~1.2]{Cardin_Viterbo_2008_Commuting_Hamiltonians_and_Hamilton-Jacobi_multi-time_equations}. Quantitative lower bounds for pairs of functions were also obtained in \cite[Corollary~1.5]{Entov_Polterovich_Zapolsky_2007_quasi-morphisms_and_PB}, in terms of the non-additivity of a symplectic quasi-state. In contrast, Corollary~\ref{cor:width-liminf} applies when the limiting map $F$ is merely continuous and the target dimension $k$ is arbitrary, and its lower bound is expressed in terms of the non-displaceability width of $F$.

\subsection{Obstructions to the quantitative thick-fiber property}
\label{ss: obst_to_BTFP}

We now present several families of closed symplectic manifolds satisfying the hypotheses of Theorem~\ref{thm: obstruction}, and hence not having the quantitative thick-fiber property. Throughout this subsection, we use the conventions
\[
\omega(X_H,\cdot)=-dH,\qquad \{F,G\}=dF(X_G).
\]
For smooth maps $u,v\colon M\to S^1$, their Poisson bracket is
\[
\{u,v\}=\Pi(u^*\theta,v^*\theta),
\]
where $\Pi$ is the Poisson bivector associated with $\omega$ and $\theta$ is the standard angular one-form on $S^1=\R/\Z$. With these conventions, the Poisson bivector associated with $dx\wedge dy$ is
\[
-\frac{\partial}{\partial x}\wedge\frac{\partial}{\partial y},
\]
and hence $\{x,y\}=-1$.

\begin{exam}
Let $\alpha>0$ be irrational, and equip $\T^4=(\R/\Z)^4$ with the symplectic form
\[\omega_\alpha=dx_1\wedge dy_1+\alpha\,dx_2\wedge dy_2.\]
The coordinate one-forms $dx_1,dy_1,dx_2,dy_2$ and the corresponding coordinate vector fields on $\R^4$ descend to $\T^4$. Define smooth maps $a,b,c\colon\T^4\to S^1=\R/\Z$ by
\[a([x_1,y_1,x_2,y_2])=x_1+\Z,\qquad b([x_1,y_1,x_2,y_2])=x_2+\Z,\qquad c([x_1,y_1,x_2,y_2])=y_1+y_2+\Z,\]
for every $[x_1,y_1,x_2,y_2]\in\T^4=\R^4/\Z^4$. The Poisson bivector associated with $\omega_\alpha$ is
\[
\Pi_\alpha
=
-\frac{\partial}{\partial x_1}\wedge\frac{\partial}{\partial y_1}
-\frac{1}{\alpha}
\frac{\partial}{\partial x_2}\wedge\frac{\partial}{\partial y_2}.
\]
Moreover,
\[
a^*\theta=dx_1,\qquad
b^*\theta=dx_2,\qquad
c^*\theta=dy_1+dy_2.
\]
Therefore, by the definition of the Poisson bracket of circle-valued functions,
\[
\{a,c\}
=
\Pi_\alpha(dx_1,dy_1+dy_2)
=-1,
\qquad
\{b,c\}
=
\Pi_\alpha(dx_2,dy_1+dy_2)
=-\frac{1}{\alpha}.
\]
Consequently,
\[
\frac{\{a,c\}}{\{b,c\}}=\alpha\notin\Q.
\]
Hence, by Theorem~\ref{thm: obstruction}, $(\T^4,\omega_\alpha)$ does not have the quantitative thick-fiber property.
\end{exam}
A generalization of the last example is the following:
\begin{exam}\label{exam: Kodaira_Thurston}
For every integer $q\geq0$, let $M_q$ be the quotient of $\R^4$, with coordinates $(x,y,z,t)$, by the equivalence relation generated by
\[\begin{aligned}
(x,y,z,t)&\sim(x+1,y,z+qy,t),\\
(x,y,z,t)&\sim(x,y+1,z,t),\\
(x,y,z,t)&\sim(x,y,z+1,t),\\
(x,y,z,t)&\sim(x,y,z,t+1).
\end{aligned}\]
Denote the quotient map by $\pi_q\colon\R^4\to M_q$. Then $M_0=\T^4$, while $M_1$ is the Kodaira--Thurston manifold; see \cite[Definition~1.1]{Evans_Kedra_2015}.

The one-form $\eta_q=dz-qx\,dy$ is invariant under these identifications and therefore descends to a smooth one-form on $M_q$. Let $\alpha>0$ be irrational, and define
\begin{equation}\label{eq: Omega_Kodaira_Thurston}
    \Omega_{\alpha,q}=\left(dx+(\alpha-1)dy\right)\wedge\eta_q+\alpha\,dy\wedge dt.
\end{equation}
Since $d\eta_q=-q\,dx\wedge dy$, we have $d\Omega_{\alpha,q}=0$. Moreover,
\[\frac{1}{2}\Omega_{\alpha,q}^2=\alpha\,dx\wedge\eta_q\wedge dy\wedge dt\neq0,\]
and hence $\Omega_{\alpha,q}$ is nondegenerate. Therefore, $\Omega_{\alpha,q}$ is a symplectic form on $M_q$.

Define smooth maps $a,b,c\colon M_q\to S^1$ by
\[a([x,y,z,t])=y-x+\Z,\qquad b([x,y,z,t])=y+\Z,\qquad c([x,y,z,t])=t+\Z,\]
for every $[x,y,z,t]\in M_q$. These maps are well defined under all the identifications defining $M_q$.

Put
\[e^1=dx+(\alpha-1)dy,\qquad e^2=\eta_q,\qquad e^3=\alpha\,dy,\qquad e^4=dt.\]
Then
\[\Omega_{\alpha,q}=e^1\wedge e^2+e^3\wedge e^4,\]
while
\[a^*\theta=-e^1+e^3,\qquad b^*\theta=\frac{1}{\alpha}e^3,\qquad c^*\theta=e^4.\]
Let $E_1,E_2,E_3,E_4$ be the frame dual to $e^1,e^2,e^3,e^4$. By our sign conventions, the Poisson bivector associated with $\Omega_{\alpha,q}$ is
\[\Pi_{\alpha,q}=-E_1\wedge E_2-E_3\wedge E_4.\]
It follows that
\[\{a,c\}=-1,\qquad \{b,c\}=-\frac{1}{\alpha}.\]
Consequently,
\[\frac{\{a,c\}}{\{b,c\}}=\alpha\notin\Q.\]
Hence, by Theorem~\ref{thm: obstruction}, $(M_q,\Omega_{\alpha,q})$ does not have the quantitative thick-fiber property for every integer $q\geq0$ and every irrational $\alpha>0$.
\end{exam}

\begin{exam}
Let $q\geq0$ be an integer, let $(N,\sigma)$ be a closed symplectic manifold, and let $\psi_1,\ldots,\psi_4\colon N\to N$ be pairwise commuting symplectomorphisms satisfying $\psi_3^q=\id_N$. We use these symplectomorphisms as monodromies to construct a fiber bundle with fiber $N$ over the manifold $M_q$ from Example~\ref{exam: Kodaira_Thurston}.

We first explain the role of the condition $\psi_3^q=\id_N$. Denote by $\tau_1,\ldots,\tau_4$ the four transformations of $\R^4$ defining $M_q$, in the order in which they appear in Example~\ref{exam: Kodaira_Thurston}, and use the convention $\tau_1\tau_2=\tau_1\circ\tau_2$. We have
\[
\tau_1\tau_2(x,y,z,t)=(x+1,y+1,z+qy+q,t)=\tau_3^q\tau_2\tau_1(x,y,z,t).
\]
Thus, $\tau_1\tau_2=\tau_3^q\tau_2\tau_1$. All pairs of the transformations $\tau_i$, except for the pair $\{\tau_1,\tau_2\}$, commute. Since the $\psi_i$ pairwise commute and $\psi_3^q=\id_N$, we have
\[
\psi_1\psi_2
=\psi_2\psi_1
=\psi_3^q\psi_2\psi_1.
\]
Thus, the assignment $\tau_i\mapsto\psi_i$ respects all the defining relations among the transformations used to construct $M_q$.

Let $M$ be the quotient of $\R^4\times N$, with coordinates $(x,y,z,t,u)$, by the equivalence relation generated by
\[
\begin{aligned}
(x,y,z,t,u)&\sim(x+1,y,z+qy,t,\psi_1(u)),\\
(x,y,z,t,u)&\sim(x,y+1,z,t,\psi_2(u)),\\
(x,y,z,t,u)&\sim(x,y,z+1,t,\psi_3(u)),\\
(x,y,z,t,u)&\sim(x,y,z,t+1,\psi_4(u)).
\end{aligned}
\]
The compatibility verified above ensures that the projection onto the first four coordinates induces a smooth fiber bundle
\[
p\colon M\to M_q
\]
with fiber $N$. Since $M_q$ and $N$ are closed, $M$ is a closed smooth manifold.

Let $\pi_q\colon\R^4\to M_q$ be the quotient map, and let $\pr_{\R^4}\colon\R^4\times N\to\R^4$ and $\pr_N\colon\R^4\times N\to N$ be the projections. Let $\alpha>0$ be irrational, let $\Omega_{\alpha,q}$ be the symplectic form on $M_q$, given in \eqref{eq: Omega_Kodaira_Thurston}, and define
\[
\widetilde{\Omega}
=
\pr_{\R^4}^*(\pi_q^*\Omega_{\alpha,q})
+\pr_N^*\sigma
\]
on $\R^4\times N$. For every $i\in\{1,\ldots,4\}$, we have $\pi_q\circ\tau_i=\pi_q$, and hence $\tau_i^*(\pi_q^*\Omega_{\alpha,q})=\pi_q^*\Omega_{\alpha,q}$. Moreover, $\psi_i^*\sigma=\sigma$ because $\psi_i$ is a symplectomorphism. Therefore, each of the four transformations defining $M$ preserves $\widetilde{\Omega}$, so $\widetilde{\Omega}$ descends to a smooth two-form $\Omega$ on $M$. The form $\widetilde{\Omega}$ is the product symplectic form determined by $\pi_q^*\Omega_{\alpha,q}$ and $\sigma$, and is therefore symplectic. Consequently, the descended form $\Omega$ is symplectic as well.

Define smooth maps $\overline{a},\overline{b},\overline{c}\colon M_q\to S^1$ by
\[
\overline{a}([x,y,z,t])=y-x+\Z,\qquad
\overline{b}([x,y,z,t])=y+\Z,\qquad
\overline{c}([x,y,z,t])=t+\Z
\]
for every $[x,y,z,t]\in M_q$. These are the circle-valued maps considered in Example~\ref{exam: Kodaira_Thurston}. Define smooth maps $a,b,c\colon M\to S^1$ by $a=\overline{a}\circ p$, $b=\overline{b}\circ p$, and $c=\overline{c}\circ p$.

Since $\Omega$ is induced by the product of $\Omega_{\alpha,q}$ and $\sigma$, the Poisson brackets of the maps pulled back from $M_q$ agree with their Poisson brackets on the base. Therefore,
\[
\{a,c\}_\Omega
=
\{\overline{a},\overline{c}\}_{\Omega_{\alpha,q}}\circ p
=-1,
\qquad
\{b,c\}_\Omega
=
\{\overline{b},\overline{c}\}_{\Omega_{\alpha,q}}\circ p
=-\frac{1}{\alpha}.
\]
Consequently,
\[
\frac{\{a,c\}_\Omega}{\{b,c\}_\Omega}
=\alpha\notin\Q.
\]
Hence, by Theorem~\ref{thm: obstruction}, $(M,\Omega)$ does not have the quantitative thick-fiber property.

For a concrete nontrivial example, take $q=0$ and let $N=\T^2$ with its standard area form. Define the symplectomorphism $\psi\colon\T^2\to\T^2$ by
\[
\psi(u)=
\begin{pmatrix}
2&1\\
1&1
\end{pmatrix}
u\pmod{\Z^2}
\]
for every $u\in\T^2$. Take $\psi_1=\psi$ and $\psi_2=\psi_3=\psi_4=\id_{\T^2}$. Let
\[
Y_\psi
=
(\R\times\T^2)/
\bigl((x,u)\sim(x+1,\psi(u))\bigr)
\]
be the mapping torus of $\psi$. The first identification defining $M$ produces $Y_\psi$ from the coordinates $(x,u)$, while the remaining three identifications take the coordinates $y$, $z$, and $t$ modulo $\Z$. Therefore,
\[
M\cong Y_\psi\times\T^3.
\]
Moreover, note that $M$ is not diffeomorphic to $\T^6$. Indeed, by \cite[Proposition~15.25 and the example on p.~479]{Friedl_Algebraic_Topology},
\[
\pi_1(M)\cong\left(\Z^2\rtimes_{\psi_*}\Z\right)\times\Z^3.
\]
This group is nonabelian, whereas $\pi_1(\T^6)\cong\Z^6$ is abelian. Therefore, $M$ is not diffeomorphic to $\T^6$.
\end{exam}

\subsection{Discussion}\label{ss: discussion}

\subsubsection{Noncompact symplectic manifolds}

It is natural to ask whether semi-symplectic quasi-states and the quantitative thick-fiber theorem can be extended to noncompact symplectic manifolds. In this setting, the choice of the space of functions requires some care. Following \cite[Section~1.2]{Lanzat_2013_quasi_morphisms_and_symplectic_quasi_states}, let $C_{cc}(M)$ denote the space of continuous functions on $M$ that are constant outside a compact subset. This space contains the constant functions and is therefore suitable for a normalization axiom. One could define a semi-symplectic quasi-state on a noncompact symplectic manifold to be a functional $\zeta\colon C_{cc}(M)\to\R$ satisfying normalization, monotonicity, the vanishing property for compactly supported functions with displaceable support, and the following one-sided Poisson-bracket inequality: there exists a constant $C>0$ such that, for every $F,G\in C_{cc}^\infty(M)$, we have
\[
\zeta(F+G)-\zeta(F)-\zeta(G)
\leq C\sqrt{\|\{F,G\}\|_\infty}.
\]

Lanzat constructed partial symplectic quasi-states on certain open convex symplectic manifolds; see \cite[Corollary~3.15]{Lanzat_2013_quasi_morphisms_and_symplectic_quasi_states}. The examples following that corollary include cotangent bundles of closed connected manifolds and open Stein manifolds; see \cite[p.~5354]{Lanzat_2013_quasi_morphisms_and_symplectic_quasi_states}. He also constructed genuine symplectic quasi-states on suitable one-point blowups of cotangent bundles; see \cite[Theorem~1.6 and Remark~6]{Lanzat_2013_quasi_morphisms_and_symplectic_quasi_states}.

On the other hand, Zapolsky proved that every stable homogeneous quasi-morphism on the universal cover of the group of compactly supported Hamiltonian diffeomorphisms of a connected complete Liouville manifold vanishes identically; see \cite[Theorem~1.2]{Zapolsky_2026_Hofer_continuous_quasi_morphisms}. Since the standard construction from stable homogeneous quasi-morphisms produces genuine symplectic quasi-states satisfying the Poisson-bracket inequality, this construction cannot produce a nontrivial genuine symplectic quasi-state on such a manifold; see also \cite[Proposition~2.8]{Zapolsky_2026_Hofer_continuous_quasi_morphisms}. This does not rule out partial or semi-symplectic quasi-states. It also does not contradict Lanzat's results: his quasi-states on the weakly exact examples above are partial, while his nontrivial genuine examples are defined on nonexact blowups.

We next ask whether Theorem~\ref{thm:thick-fiber} admits an analogue for noncompact symplectic manifolds. Let $\zeta\colon C_{cc}(M)\to\R$ be a functional of the type described above. Before discussing such an extension, we note that the notions of $\zeta$-pseudoheavy and $\zeta$-superheavy subsets can be adapted verbatim to the noncompact setting; the notion of a superheavy subset was adapted by Lanzat in \cite[Definition~4.6]{Lanzat_2013_quasi_morphisms_and_symplectic_quasi_states}. Displaceability should be understood with respect to compactly supported Hamiltonian diffeomorphisms. As in the closed case, every $\zeta$-pseudoheavy subset intersects every compact $\zeta$-superheavy subset, and the proof is the same as that of \cite[Proposition~1.6]{Kawasaki_Orita_2021_seudoheavy}.

A plausible noncompact analogue of Theorem~\ref{thm:thick-fiber} can now be formulated for smooth maps $F\colon M\to\R^k$ whose coordinate functions belong to $C_{cc}^\infty(M)$. Such a map is constant outside a compact subset; denote its value there by $p_\infty\in\R^k$. It seems plausible that suitable modifications of the proof of Theorem~\ref{thm:thick-fiber} would produce a $\zeta$-pseudoheavy closed thick fiber
\[
F^{-1}\left(\overline{B}_{y_0}(R_0)\right),
\qquad
R_0=4Ck^{3/2}\sqrt{pb(F)}.
\]
This conclusion may be trivial if $p_\infty\in\overline{B}_{y_0}(R_0)$, since the thick fiber then contains the complement of a compact subset and is automatically non-displaceable by compactly supported Hamiltonian diffeomorphisms. One possible way to exclude this case is to assume that there exists a compact $\zeta$-superheavy subset $K\subseteq M$ satisfying
\[
\dist\left(p_\infty,F(K)\right)>2R_0.
\]
The pseudoheavy thick fiber would then intersect $K$. Consider a point
\[
x\in K\cap F^{-1}\left(\overline{B}_{y_0}(R_0)\right).
\]
We have
\[
\|y_0-p_\infty\|
\geq\|F(x)-p_\infty\|-\|F(x)-y_0\|
>2R_0-R_0=R_0.
\]
Consequently, $p_\infty\notin\overline{B}_{y_0}(R_0)$. Since $F=p_\infty$ outside a compact subset, it follows that
\[
F^{-1}\left(\overline{B}_{y_0}(R_0)\right)
\]
is compact. This suggests a meaningful noncompact version of the quantitative thick-fiber theorem: a compact superheavy subset whose image is quantitatively separated from the value of the map at infinity should force the existence of a compact pseudoheavy, and hence non-displaceable, thick fiber.

\subsubsection{Riemannian manifolds as target spaces}

The definition of a thick fiber extends naturally to maps with values in a Riemannian manifold. Let $(N,g)$ be a $k$-dimensional Riemannian manifold, and let $F\colon M\to N$ be a smooth map. For every $y\in N$ and every $r\geq0$, one can consider the thick fiber
\[F^{-1}\left(\overline{B}_y^g(r)\right),\]
where $\overline{B}_y^g(r)$ is the closed metric ball of radius $r$ centered at $y$. To formulate a quantitative thick-fiber problem in this setting, we first need a coordinate-independent measure of the failure of the components of $F$ to Poisson commute.

Let $\Pi$ be the Poisson bivector associated with $\omega$. Define the section $\cP_F$ of the pullback bundle $F^*\left(\bigwedge\nolimits^2 TN\right)\to M$ as the pushforward of $\Pi$ under $F$. More precisely, for every $x\in M$ and every $\alpha,\beta\in T_{F(x)}^*N$, define
\[
(\cP_F)_x(\alpha,\beta)
=
\Pi_x\left((d_xF)^*\alpha,(d_xF)^*\beta\right).
\]
The Riemannian metric $g$ induces a norm on $\bigwedge^2 TN$. Define the \textbf{Poisson-bracket size of $F$ with respect to $g$} by
\[
pb_g(F)=\max_{x\in M}\|(\cP_F)_x\|_g.
\]

In local coordinates $(y_1,\ldots,y_k)$ on $N$, we have
\[\mathcal{P}_F(x)=\sum_{1\leq i<j\leq k}\{y_i\circ F,y_j\circ F\}(x)\left.\frac{\partial}{\partial y_i}\wedge\frac{\partial}{\partial y_j}\right|_{F(x)}.\]
Thus, $pb_g(F)$ measures simultaneously the Poisson brackets of all the local coordinate functions of $F$, in a coordinate-independent way.

When $N=\R^k$ is equipped with the Euclidean metric and $F=(F_1,\ldots,F_k)$, we have
\[pb_g(F)=\max_{x\in M}\left(\sum_{1\leq i<j\leq k}|\{F_i,F_j\}(x)|^2\right)^{1/2}.\]
Consequently,
\[pb(F)\leq pb_g(F)\leq\sqrt{\frac{k(k-1)}{2}}\,pb(F).\]
In particular, $pb_g(F)=0$ if and only if the coordinate functions of $F$ pairwise Poisson commute.

Theorem~\ref{thm:thick-fiber} also gives a local result for Riemannian targets. Suppose that $F(M)$ is contained in a coordinate chart whose distortion in both directions is bounded by $L\geq1$. Applying Theorem~\ref{thm:thick-fiber} in this chart, we obtain a point $y\in F(M)$ such that
\[
F^{-1}\left(\overline{B}_y^g\left(4Ck^{3/2}L^2\sqrt{pb_g(F)}\right)\right)
\]
contains a $\zeta$-pseudoheavy subset and is therefore non-displaceable. Indeed, the chart increases the Poisson-bracket size by at most a factor $L^2$, which contributes a factor $L$ after taking the square root, while returning to the Riemannian distance contributes another factor $L$.

For a general Riemannian target, the local coordinate argument does not immediately give a global estimate. Such an estimate would require a quantitatively controlled collection of coordinate neighborhoods covering the relevant part of $N$, together with a way to combine the corresponding local estimates. The resulting constant should depend on the dimension and geometric quantities such as the injectivity radius, curvature bounds, or the bi-Lipschitz distortion of the coordinate charts. Since $F(M)$ is compact, only the geometry of a neighborhood of $F(M)$ should be relevant when $N$ is noncompact.

The square-root dependence remains natural in this setting. Indeed, if the metric is rescaled by $g_\lambda=\lambda^2g$ for some $\lambda>0$, then
\[
pb_{g_\lambda}(F)=\lambda^2pb_g(F),
\]
while
\[
d_{g_\lambda}=\lambda d_g.
\]
Therefore, a radius proportional to $\sqrt{pb_g(F)}$ transforms in the same way as the Riemannian distance. This suggests that the square-root dependence on the Poisson-bracket size should remain unchanged, while the dependence on the dimension should be supplemented by a quantitative constant determined by the geometry of the target.

It would be interesting to determine which geometric assumptions on $(N,g)$ guarantee the existence of a constant $C(N,g)>0$ such that, for every smooth map $F\colon M\to N$, there exists $y\in F(M)$ for which
\[
F^{-1}\left(\overline{B}_y^g\left(C(N,g)\sqrt{pb_g(F)}\right)\right)
\]
is non-displaceable. In particular, one may ask whether manifolds of bounded geometry admit an estimate depending only on the dimension and the relevant bounded-geometry constants, and whether the topology of the target produces additional obstructions beyond those already present for Euclidean targets.

\subsubsection{The displacement energy of fibers}
\label{sss:fiber-displacement-energy}

Let $X\subseteq M$ be a closed subset. The \textbf{displacement energy} of $X$ is defined by
\[
e(X)
=
\inf\left\{
\|H\|_{\mathrm{Hofer}}:
H\in C^\infty(S^1\times M),\
\varphi_H^1(X)\cap X=\varnothing
\right\},
\]
where $\varphi_H^1$ is the time-one map of the Hamiltonian flow generated by $H$, and
\[
\|H\|_{\mathrm{Hofer}}
=
\int_{S^1}
\left(
\max_{x\in M}H(t,x)-\min_{x\in M}H(t,x)
\right)\,dt.
\]
We use the convention that $e(X)=\infty$ if $X$ is non-displaceable.

For every continuous map $F\colon M\to\R^k$, define the maximal displacement energy of its fibers by
\[
e_{\mathrm{fib}}(F)
=
\sup_{p\in\R^k}e\left(F^{-1}(p)\right).
\]
This quantity cannot be usefully estimated in terms of $pb(F)$ alone. Indeed, if $F$ is a smooth embedding and $M$ has positive dimension, then every nonempty fiber is a singleton, and hence
\[
e_{\mathrm{fib}}(F)=0.
\]
Moreover, for every $\lambda>0$, the maps $F$ and $\lambda F$ have the same fibers, up to relabeling, whereas
\[
pb(\lambda F)=\lambda^2pb(F).
\]
On the other hand, if $pb(F)=0$, then Theorem~\ref{thm:EP-fiber} gives a non-displaceable fiber, and hence
\[
e_{\mathrm{fib}}(F)=\infty.
\]
Thus, an estimate for $e_{\mathrm{fib}}(F)$ requires additional assumptions on the map or on the geometry of its fibers.

The displacement energies of thick fibers can be recorded by the function $\mathcal E_F\colon[0,\infty)\to[0,\infty]$ defined by
\[
\mathcal E_F(r)
=
\sup_{p\in\R^k}
e\left(
F^{-1}\left(\overline{B}_p(r)\right)
\right)
\]
for every $r\geq0$. The function $\mathcal E_F$ is nondecreasing, and
\[
\mathcal E_F(0)=e_{\mathrm{fib}}(F).
\]
Moreover, using the compactness of $M$, one can verify that the non-displaceability width is given by
\[
w(F)=\inf\left\{r>0:\mathcal E_F(r)=\infty\right\}.
\]

In this notation, Theorem~\ref{thm:thick-fiber} gives
\[
\mathcal E_F\left(
4Ck^{3/2}\sqrt{pb(F)}
\right)
=
\infty.
\]
The displacement-energy profile also satisfies a $C^0$-stability property. If $F,G\colon M\to\R^k$ are continuous and
\[
\|F-G\|_\infty\leq\epsilon,
\]
then, for every $r\geq0$,
\[
\mathcal E_F(r)\leq\mathcal E_G(r+\epsilon),
\qquad
\mathcal E_G(r)\leq\mathcal E_F(r+\epsilon).
\]
This follows from the corresponding inclusions between the thick fibers and the monotonicity of displacement energy.

\begin{question}
Under which additional assumptions on a smooth map $F\colon M\to\R^k$ can one estimate $e_{\mathrm{fib}}(F)$? More generally, can one estimate $\mathcal E_F(r)$ for $0\leq r<w(F)$ using symplectic information about the fibers or the subsets containing them?
\end{question}

\paragraph*{Acknowledgements.}
I am grateful to Yaniv Ganor for many interesting mathematical discussions and for his advice concerning the development of this paper and the use of AI in the research process. I am also grateful to Jun Zhang for sharing preliminary versions of his work ``Contact Rigidity and Poisson Bracket Invariants,'' for helpful discussions concerning the relationship between our results, and for the open and collegial manner in which we coordinated the presentation of our related work. I would also like to thank Michael Entov, Leonid Polterovich, Egor Shelukhin and David Keren Yaar for useful discussions.

I was partially supported by the Fondation Courtois during my postdoctoral fellowship at the Centre de Recherches Math\'ematiques, Universit\'e de Montr\'eal, and by the Jenny Panitch Beckow Memorial Scholarship, Canada, 2026. I also thank the Centre de Recherches Math\'ematiques for its warm hospitality and excellent working environment.

\paragraph*{AI disclosure.}
The author conceived and directed the research program of this paper. He formulated the quantitative thick-fiber property and the notion of a semi-symplectic quasi-state, identified the one-sided Poisson-bracket inequality as the relevant hypothesis, and directed the study of the existence, obstruction, optimality, and approximation questions considered here. He identified the general mechanism underlying Theorem~\ref{thm: obstruction} and developed the final $(k+1)$-color partition-of-unity construction in Subsection~\ref{ss:partition}, inspired by \cite[Theorem~6.1]{DGPZ_2024_Symp_top_and_IVMs} and Milnor's lemma; see \cite[Lemma~2.4]{Palais_homotopy_theory_of_infinite_dim_mfds}.

OpenAI's ChatGPT and Codex contributed preliminary mathematical arguments and examples. In particular, they produced an initial partition-of-unity argument using $2^k$ colors for the proof of Theorem~\ref{thm:thick-fiber} and the initial irrational-torus example that led to Theorem~\ref{thm: obstruction}. They also assisted with preliminary arguments for Theorem~\ref{thm:lower-bound-thick-fiber-radius}, the applications of the non-displaceability width in Section~\ref{ss: non-disp width}, and the examples in Section~\ref{ss: obst_to_BTFP}. They were further used for checking computations, locating references, drafting preliminary proofs, editing the English, and preparing the \LaTeX{} source.

The author critically assessed the AI-generated material, selected the parts to be retained, substantially reworked and integrated them into his research program, corrected errors, and independently verified all mathematical arguments and references. He determined the final statements, proofs, organization, and exposition of the paper.

\section{Proofs}\label{s: proofs}

\subsection{Proof of Proposition~\ref{prop:ssqs_on_surfaces}}\label{ss: proof of prop: ssqs_on_surfaces}

The proof of Proposition~\ref{prop:ssqs_on_surfaces} is based on two facts. The first is that partial symplectic quasi-states are $1$-Lipschitz with respect to the uniform norm. The second is the following result due to Zapolsky \cite{Zapolsky_2010_almost_commuting_inD2}.

\begin{thm}[{Zapolsky, \cite[Theorem~1.3]{Zapolsky_2010_almost_commuting_inD2}}]\label{thm:Zap-approx}
Let $F,G\in C^\infty(\Sigma)$ and let $\epsilon\geq0$. If
\[\int_\Sigma|\{F,G\}|\sigma\leq\epsilon,\]
then there exist smooth functions $F',G'\in C^\infty(\Sigma)$ satisfying $\{F',G'\}=0$ and
\[\|F-F'\|_\infty\leq\sqrt{\frac{\epsilon}{2}},\qquad \|G-G'\|_\infty\leq\sqrt{\frac{\epsilon}{2}}.\]
\end{thm}

\begin{proof}[Proof of Proposition~\ref{prop:ssqs_on_surfaces}]
Since $\zeta$ is a partial symplectic quasi-state, it satisfies the normalization, monotonicity, and vanishing properties. Let us show that $\zeta$ satisfies the one-sided Poisson-bracket inequality with constant $2\sqrt{2A}$.

Let $F,G\in C^\infty(\Sigma)$. We have
\[\int_\Sigma|\{F,G\}|\sigma\leq\int_\Sigma\|\{F,G\}\|_\infty\sigma=A\|\{F,G\}\|_\infty.\]
By Theorem~\ref{thm:Zap-approx}, there exist Poisson-commuting smooth functions $F',G'\in C^\infty(\Sigma)$ such that
\[\|F-F'\|_\infty\leq\sqrt{\frac{A}{2}}\sqrt{\|\{F,G\}\|_\infty},\qquad \|G-G'\|_\infty\leq\sqrt{\frac{A}{2}}\sqrt{\|\{F,G\}\|_\infty}.\]
In particular,
\[\|(F+G)-(F'+G')\|_\infty\leq\sqrt{2A}\sqrt{\|\{F,G\}\|_\infty}.\]
Since $F'$ and $G'$ Poisson commute, the quasi-subadditivity of $\zeta$ implies that
\[\zeta(F'+G')\leq\zeta(F')+\zeta(G').\]
Finally, since $\zeta$ is $1$-Lipschitz with respect to the uniform norm, we obtain
\begin{align*}
\zeta(F+G)&\leq\zeta(F'+G')+\sqrt{2A}\sqrt{\|\{F,G\}\|_\infty}\\
&\leq\zeta(F')+\zeta(G')+\sqrt{2A}\sqrt{\|\{F,G\}\|_\infty}\\
&\leq\zeta(F)+\zeta(G)+2\sqrt{2A}\sqrt{\|\{F,G\}\|_\infty}.
\end{align*}
Therefore,
\[\zeta(F+G)-\zeta(F)-\zeta(G)\leq2\sqrt{2A}\sqrt{\|\{F,G\}\|_\infty},\]
as required.
\end{proof}

\subsection{Proof of Theorem~\ref{thm:lower-bound-thick-fiber-radius}}\label{ss: Proof of the lower-bound-thick-fiber-radius}

The first step in the proof of Theorem~\ref{thm:lower-bound-thick-fiber-radius} is to establish the nontriviality of the non-displaceability width.

\begin{lemma}\label{lemma:partition-width}
Let $U_1,\ldots,U_k$ be a finite open cover of $M$ by displaceable sets, and let $\rho_1,\ldots,\rho_k\colon M\to[0,1]$ be a smooth partition of unity subordinate to this cover. Define the smooth map $F\colon M\to\R^k$ by $F(x)=\bigl(\rho_1(x),\ldots,\rho_k(x)\bigr)$ for every $x\in M$. Then $w(F)\geq\frac{1}{2k}$.
\end{lemma}
\begin{proof}
We prove that, for every $p\in\R^k$ and every $r\in(0,1/(2k))$, the set $F^{-1}\left(\overline{B}_p(r)\right)$ is displaceable. Fix $p\in\R^k$ and $r\in(0,1/(2k))$, and put $X=F^{-1}\left(\overline{B}_p(r)\right)$. If $X$ is empty, then it is displaceable. Suppose that $X$ is nonempty, and choose $x_0\in X$. Since $\sum_{i=1}^k\rho_i(x_0)=1$, there exists $i\in\{1,\ldots,k\}$ such that $\rho_i(x_0)\geq\frac{1}{k}$. For every $x\in X$, we have
\[
|\rho_i(x)-\rho_i(x_0)|\leq\|F(x)-F(x_0)\|\leq\|F(x)-p\|+\|p-F(x_0)\|\leq2r<\frac{1}{k}.
\]
It follows that $\rho_i(x)\geq\rho_i(x_0)-2r\geq\frac{1}{k}-2r>0$. This shows that
\[
X\subseteq\{\rho_i>0\}\subseteq\supp\rho_i\subseteq U_i.
\]
Since $U_i$ is displaceable, the set $X$ is displaceable. Consequently, $w(F)\geq\frac{1}{2k}$.
\end{proof}

\begin{rem}
By combining Theorem~\ref{thm:thick-fiber} with Lemma~\ref{lemma:partition-width}, we obtain a lower bound for the Poisson-bracket invariant of open covers by displaceable sets; see \cite{Entov_Polterovich_Zapolsky_2007_quasi-morphisms_and_PB}. However, this lower bound does not improve the lower bound obtained in that paper.
\end{rem}
\begin{proof}[Proof of Theorem~\ref{thm:lower-bound-thick-fiber-radius}]
Every point of $M$ has a displaceable open neighborhood contained in a Darboux chart. By the compactness of $M$, there exist an integer $q\geq2$ and a finite open cover $U_1,\ldots,U_q$ of $M$ by displaceable sets. Let $\rho_1,\ldots,\rho_q\colon M\to[0,1]$ be a smooth partition of unity subordinate to this cover. Define the smooth map $\Phi\colon M\to\R^q$ by
\[
\Phi(x)=\bigl(\rho_1(x),\ldots,\rho_q(x)\bigr)
\]
for every $x\in M$, and put $\delta=pb(\Phi)$. By Lemma~\ref{lemma:partition-width}, we have $w(\Phi)\geq\frac{1}{2q}>0$. Since $(M,\omega)$ has the quantitative thick-fiber property with respect to $f$ and $g$, we also have
\[
w(\Phi)\leq f(q)g(\delta).
\]
If $\delta=0$, then $g(\delta)=g(0)=0$, and hence the preceding inequality gives $w(\Phi)\leq0$, contradicting the fact that $w(\Phi)>0$. Therefore, $\delta>0$. Moreover, since $w(\Phi)>0$ and $f(q)>0$, the same inequality implies that $g(\delta)>0$. We prove the required bounds in two stages.
\begin{enumerate}
\item Let us first prove the lower bound for $g$. Consider $c_g=\frac{w(\Phi)}{f(q)\sqrt{\delta}}>0$. Let $t\geq0$. If $t=0$, then $g(t)=0=c_g\sqrt{t}$. Suppose that $t>0$, and put $\lambda=\sqrt{\frac{t}{\delta}}>0$. We have $pb(\lambda\Phi)=\lambda^2pb(\Phi)=t$, and Proposition~\ref{prop:width-stability} gives $w(\lambda\Phi)=\lambda w(\Phi)$. Applying the quantitative thick-fiber property to the smooth map $\lambda\Phi\colon M\to\R^q$, we obtain
\[
\lambda w(\Phi)=w(\lambda\Phi)\leq f(q)g\bigl(pb(\lambda\Phi)\bigr)=f(q)g(t).
\]
Therefore,
\[
g(t)\geq\frac{\lambda w(\Phi)}{f(q)}
=\frac{w(\Phi)}{f(q)\sqrt{\delta}}\sqrt{t}
=c_g\sqrt{t}.
\]
This proves the required lower bound for $g$.

\item Let us now prove the lower bound for $f$. Let $k\geq q$, and put $m=\lfloor k/q\rfloor$. Define the smooth map $F_k\colon M\to\R^k$ by
\[
F_k(x)=\left(\underbrace{\Phi(x),\ldots,\Phi(x)}_{m\text{ copies}},0,\ldots,0\right)
\]
for every $x\in M$. We have $pb(F_k)=\delta$. By the quantitative thick-fiber property, there exists $y\in F_k(M)$ such that
\[
F_k^{-1}\left(\overline{B}_y\left(f(k)g(\delta)\right)\right)
\]
is non-displaceable. Choose $x_0\in M$ such that $y=F_k(x_0)$. For every $x\in M$, we have
\[
\|F_k(x)-F_k(x_0)\|=\sqrt{m}\|\Phi(x)-\Phi(x_0)\|.
\]
Therefore,
\[
F_k^{-1}\left(\overline{B}_y\left(f(k)g(\delta)\right)\right)
=
\Phi^{-1}\left(
\overline{B}_{\Phi(x_0)}
\left(\frac{f(k)g(\delta)}{\sqrt{m}}\right)
\right).
\]
Since this set is non-displaceable, the definition of the non-displaceability width and Lemma~\ref{lemma:partition-width} give
\[
\frac{1}{2q}\leq w(\Phi)\leq\frac{f(k)g(\delta)}{\sqrt{m}}.
\]
Since $m=\lfloor k/q\rfloor\geq k/(2q)$, we obtain
\[
f(k)\geq\frac{1}{2q\sqrt{2q}\,g(\delta)}\sqrt{k}
\]
for every integer $k\geq q$. Define
\[
c_f=
\min\left\{
\frac{1}{2q\sqrt{2q}\,g(\delta)},
\min_{1\leq k<q}\frac{f(k)}{\sqrt{k}}
\right\}>0.
\]
Then $f(k)\geq c_f\sqrt{k}$ for every integer $k\geq1$. This proves the required lower bound for $f$ and completes the proof.
\end{enumerate}
\end{proof}

\subsection{Proof of Theorem~\ref{thm: obstruction}}\label{ss: Proof of obstruction}

The proof of Theorem~\ref{thm: obstruction} uses some basic properties of the Poisson bracket for circle-valued functions. Let us recall the definition and record these properties.

We use the notation $S^1=\R/\Z$ and $\theta\in\Omega^1(S^1)$ introduced in the Introduction. Recall that $\Pi$ is the Poisson bivector associated with $\omega$, characterized by
\[\Pi(dF,dG)=\{F,G\}\]
for every $F,G\in C^\infty(M)$. If $f,g\colon M\to S^1$ are smooth maps, recall that
\[\{f,g\}=\Pi(f^*\theta,g^*\theta).\]
If $f\colon M\to S^1$ and $G\colon M\to\R$ are smooth, define
\[\{f,G\}=\Pi(f^*\theta,dG).\]
For every smooth function $h\colon S^1\to\R$, denote by $h'\colon S^1\to\R$ the smooth function determined by $dh=h'\theta$. If $f,g\colon M\to S^1$ are smooth, then
\[
\{f,h\circ g\}=h'(g)\{f,g\}.
\]
Moreover, if $h_1,h_2\colon S^1\to\R$ are smooth, then
\[
\{h_1\circ f,h_2\circ g\}=h_1'(f)h_2'(g)\{f,g\}.
\]
Since $\Pi$ is skew-symmetric, we have $\{f,f\}=0$. Consequently, $\{h_1\circ f,h_2\circ f\}=0$.

Let $m,n\in\Z$, and define the smooth map $w\colon M\to S^1$ by $w=mf+ng$, where addition and integer multiplication are taken pointwise in $S^1=\R/\Z$. Then
\[w^*\theta=mf^*\theta+ng^*\theta.\]
Consequently, for every smooth map $q\colon M\to S^1$ and every smooth function $H\colon M\to\R$, we have
\[\{w,q\}=m\{f,q\}+n\{g,q\},\qquad \{w,H\}=m\{f,H\}+n\{g,H\}.\]

Finally, let $f\colon M\to S^1$ and $H\colon M\to\R$ be smooth, let $(\varphi_H^t)_{t\in\R}$ be the Hamiltonian flow of $H$, and let $x\in M$. Define the smooth path $\gamma_x\colon\R\to S^1$ by
\[\gamma_x(t)=f(\varphi_H^t(x))\]
for every $t\in\R$. If $\widetilde{\gamma}_x\colon\R\to\R$ is a smooth lift of $\gamma_x$, then
\[\frac{d}{dt}\widetilde{\gamma}_x(t)=\{f,H\}(\varphi_H^t(x))\]
for every $t\in\R$. In particular, if $\{f,H\}=0$, then the Hamiltonian flow of $H$ preserves $f$, that is,
\[f(\varphi_H^t(x))=f(x)\]
for every $x\in M$ and every $t\in\R$.

\begin{proof}[Proof of Theorem~\ref{thm: obstruction}]
Suppose, by contradiction, that $(M,\omega)$ has the quantitative thick-fiber property with respect to functions $f\colon\N\to(0,\infty)$ and $g\colon[0,\infty)\to[0,\infty)$, where $g$ is continuous and vanishes at $0$. Since $g$ is continuous and $g(0)=0$, there exists $\eta>0$ such that, for every $t\in[0,\eta]$, we have $g(t)\leq 1/(4f(4))$.

Since $A/B$ is irrational, the pigeonhole principle implies that there exist $m,n\in\Z$ such that
\[0<mA+nB<\frac{\eta}{4\pi^2}.\]
Put $\delta=mA+nB$, and define the smooth map $w\colon M\to S^1$ by
\[w(x)=ma(x)+nb(x)\]
for every $x\in M$, where addition is taken in $S^1=\R/\Z$. By the linearity property above, we have
\[\{w,c\}=m\{a,c\}+n\{b,c\}=mA+nB=\delta.\]
Define the smooth map $F\colon M\to\R^4$ by
\[F(x)=\bigl(\cos(2\pi w(x)),\sin(2\pi w(x)),\cos(2\pi c(x)),\sin(2\pi c(x))\bigr)\]
for every $x\in M$.

Note that $pb(F)<\eta$. Indeed, the first two coordinates of $F$ Poisson commute with each other, and the last two coordinates also Poisson commute with each other. Moreover, by the chain-rule property above, we have
\[\{\cos(2\pi w),\cos(2\pi c)\}=4\pi^2\sin(2\pi w)\sin(2\pi c)\{w,c\}=4\pi^2\sin(2\pi w)\sin(2\pi c)\delta.\]
Therefore,
\[\|\{\cos(2\pi w),\cos(2\pi c)\}\|_\infty\leq4\pi^2\delta.\]
Similarly, the uniform norm of the Poisson bracket of every other pair of coordinates of $F$ is bounded from above by $4\pi^2\delta$. Hence $pb(F)\leq4\pi^2\delta<\eta$.

Since $(M,\omega)$ has the quantitative thick-fiber property with respect to $f$ and $g$, there exists $y_0\in F(M)$ such that
\[F^{-1}\left(\overline{B}_{y_0}\left(f(4)g\left(pb(F)\right)\right)\right)\]
is non-displaceable. Since $pb(F)<\eta$, we have $g(pb(F))\leq 1/(4f(4))$, and therefore $f(4)g(pb(F))\leq\frac{1}{4}$. Consequently, the set
\[K=F^{-1}\left(\overline{B}_{y_0}\left(\frac{1}{4}\right)\right)\]
is non-displaceable, since it contains the preceding non-displaceable set.

Since $y_0\in F(M)$, there exist $y_1,y_2\in S^1=\R/\Z$ such that
\[y_0=\bigl(\cos(2\pi y_1),\sin(2\pi y_1),\cos(2\pi y_2),\sin(2\pi y_2)\bigr).\]
Consider the symmetric open intervals
\[I=\left\{y_1+t\in S^1:-\frac{1}{8}<t<\frac{1}{8}\right\},\qquad J=\left\{y_2+t\in S^1:-\frac{1}{8}<t<\frac{1}{8}\right\},\]
where addition is taken modulo $\Z$.

Let us show that $K\subseteq(w,c)^{-1}(I\times J)$. Let $x\in K$. There exist $t_1,t_2\in[-1/2,1/2]$ such that $w(x)=y_1+t_1$ and $c(x)=y_2+t_2$ in $S^1=\R/\Z$. By the definition of $K$, we have $\|F(x)-y_0\|\leq\frac{1}{4}$. Using the identity
\[(\cos u-\cos v)^2+(\sin u-\sin v)^2=4\sin^2\left(\frac{u-v}{2}\right),\]
which holds for every $u,v\in\R$, we obtain
\[4\sin^2(\pi t_1)+4\sin^2(\pi t_2)=\|F(x)-y_0\|^2\leq\frac{1}{16}.\]
Consequently, $2|\sin(\pi t_i)|\leq\frac{1}{4}$ for every $i\in\{1,2\}$. Recall that, for every $u\in[0,1/2]$, we have $2u\leq\sin(\pi u)$. Since $|t_i|\leq\frac{1}{2}$, we deduce that
\[4|t_i|\leq2\sin(\pi|t_i|)=2|\sin(\pi t_i)|\leq\frac{1}{4}.\]
Therefore, $|t_i|\leq\frac{1}{16}<\frac{1}{8}$ for every $i\in\{1,2\}$. Thus $w(x)\in I$ and $c(x)\in J$. 

This shows that
\[K\subseteq(w,c)^{-1}(I\times J).\]

Choose a smooth function $h\colon S^1\to\R$ whose derivative equals $\frac{1}{2\delta}$ on a neighborhood of $\overline{J}$. Such a function exists because $\overline{J}$ is contained in a proper open arc of $S^1$. Define the Hamiltonian $H\colon M\to\R$ by $H=h\circ c$, and let $(\varphi_H^t)_{t\in\R}$ be its Hamiltonian flow. We show that its time-$1$ map $\varphi_H^1$ displaces $K$ from itself.

Let $x\in K$. Since
\[\{c,H\}=h'(c)\{c,c\}=0,\]
the Hamiltonian flow of $H$ preserves $c$. Therefore,
\[c(\varphi_H^t(x))=c(x)\in J\]
for every $t\in\R$.

Define the smooth path $\gamma_x\colon[0,1]\to S^1$ by $\gamma_x(t)=w(\varphi_H^t(x))$ for every $t\in[0,1]$, and let $\widetilde{\gamma}_x\colon[0,1]\to\R$ be a smooth lift of $\gamma_x$. By the properties above, we have
\[\frac{d}{dt}\widetilde{\gamma}_x(t)=\{w,H\}(\varphi_H^t(x))=h'(c(\varphi_H^t(x)))\{w,c\}(\varphi_H^t(x))=\frac{1}{2\delta}\delta=\frac{1}{2}\]
for every $t\in[0,1]$. Integrating this equality gives $\widetilde{\gamma}_x(1)=\widetilde{\gamma}_x(0)+\frac{1}{2}$. Passing to $S^1=\R/\Z$, we obtain
\[w(\varphi_H^1(x))=w(x)+\frac{1}{2}.\]
Consequently,
\[\varphi_H^1(K)\subseteq(w,c)^{-1}\left(\left(I+\frac{1}{2}\right)\times J\right).\]
The interval $I$ has length $\frac{1}{4}$, and hence
\[I\cap\left(I+\frac{1}{2}\right)=\varnothing.\]
It follows that
\[\varphi_H^1(K)\cap K=\varnothing.\]
Thus $K$ is displaceable, contradicting its non-displaceability. Therefore, $(M,\omega)$ does not have the quantitative thick-fiber property.
\end{proof}

\subsection{Proof of Theorem~\ref{thm:thick-fiber}}\label{ss: proof of thick-fiber}

In this section, we use several norms on $\R^k$, on spaces of functions, and on spaces of $1$-forms. We introduce them here.

We equip $\R^k$ with the coordinates $y=(y_1,\ldots,y_k)$ associated with the standard basis. For every $y=(y_1,\ldots,y_k)\in\R^k$, define
\[\|y\|_\infty=\max_i |y_i|.\]
For every $p\in[1,\infty)$, define the $\ell^p$-norm on $\R^k$ by
\[\|y\|_p=\left(\sum_{i=1}^k |y_i|^p\right)^{1/p},\]
for every $y=(y_1,\ldots,y_k)\in\R^k$. The $\ell^2$-norm is the Euclidean norm, and we denote it by $\|y\|=\|y\|_2$.

For every covector $\xi=\sum_{i=1}^k\xi_i\,dy_i$, put $\|\xi\|_1=\sum_{i=1}^k|\xi_i|$. Let $\Omega\subseteq\R^k$ be a domain. For every bounded continuous function $f\colon\Omega\to\R$, define
\[\|f\|_{\infty,\Omega}=\sup_{x\in\Omega}|f(x)|.\]
For every continuous $1$-form $\alpha$ on $\Omega$ satisfying $\sup_{x\in\Omega}\|\alpha(x)\|_1<\infty$, write $\alpha(x)=\sum_{i=1}^k\alpha_i(x)\,dy_i$ and define
\[\|\alpha\|_{1,\infty}=\sup_{x\in\Omega}\|\alpha(x)\|_1=\sup_{x\in\Omega}\sum_{i=1}^k|\alpha_i(x)|.\]
In particular, for every $u\in C^1(\Omega)$ satisfying $\sup_{x\in\Omega}\|du(x)\|_1<\infty$,
\[\|du\|_{1,\infty}=\sup_{x\in\Omega}\sum_{i=1}^k|\partial_i u(x)|.\]
When the domain is clear, we omit it from the notation for the uniform norm.

\subsubsection{Special partition of unity}
\label{ss:partition}

In this subsection, we construct a special partition of unity on $\R^k$ that will be crucial in the proof of Theorem~\ref{thm:thick-fiber}. The construction uses the Freudenthal--Kuhn triangulation of $\R^k$; see \cite{Freudenthal_1942_simplicial_decompositions,Kuhn_1960_combinatorial_lemmas}. The use of $k+1$ colors, with the functions of each color having pairwise disjoint positivity sets, was inspired by the proof of \cite[Theorem~6.1]{DGPZ_2024_Symp_top_and_IVMs} and by the covering-dimension argument in \cite{Karasev_2014_Covering_dimension_using_toric_varieties}. Both arguments use a lemma of Milnor, stated in \cite[Lemma~2.4]{Palais_homotopy_theory_of_infinite_dim_mfds}, which decomposes a suitable refinement of an open cover of a $k$-dimensional paracompact space into $k+1$ families of pairwise disjoint open sets. The explicit construction below is adapted to the quantitative differential estimates needed here.

Denote by $S_k$ the symmetric group of the set $\{1,\ldots,k\}$. Thus, the elements of $S_k$ are the permutations of $\{1,\ldots,k\}$. For every $\sigma\in S_k$, define the simplex
\[\Delta_\sigma=\left\{t=(t_1,\ldots,t_k)\in[0,1]^k:1\geq t_{\sigma(1)}\geq\cdots\geq t_{\sigma(k)}\geq0\right\}.\]
The collection $(\Delta_\sigma)_{\sigma\in S_k}$ forms a triangulation of $[0,1]^k$.

Fix $h>0$ for the remainder of this subsection. For every $m\in\Z^k$, define the affine transformation $T_m\colon\R^k\to\R^k$ by $T_m(v)=h(m+v)$ for every $v\in\R^k$. For every $m\in\Z^k$ and every $\sigma\in S_k$, define $\Delta_{m,\sigma}=T_m(\Delta_\sigma)$. The collection $(\Delta_{m,\sigma})_{m\in\Z^k,\,\sigma\in S_k}$ forms a triangulation of $\R^k$, called the \textbf{Freudenthal--Kuhn triangulation at scale $h$}. Its elements are called the \textbf{Freudenthal--Kuhn simplices at scale $h$}.

\begin{defin}
For every $v\in h\Z^k$, define $\beta_v\colon\R^k\to[0,1]$ to be the unique continuous function that is affine on every simplex in $(\Delta_{m,\sigma})_{m\in\Z^k,\,\sigma\in S_k}$ and satisfies $\beta_v(w)=1$ if $w=v$ and $\beta_v(w)=0$ otherwise, for every $w\in h\Z^k$.
\end{defin}

We will show that the family $(\beta_v)_{v\in h\Z^k}$ is a locally finite nonnegative partition of unity on $\R^k$. This collection of functions can be described explicitly using the following norm. For every $x=(x_1,\ldots,x_k)\in\R^k$, define its \textbf{Freudenthal--Kuhn norm} by
\begin{equation}\label{eq:FK-norm}
\|x\|_{\mathrm{FK}}=\max\{0,x_1,\ldots,x_k\}-\min\{0,x_1,\ldots,x_k\}.
\end{equation}
If $x_0=0$, then $\|x\|_{\mathrm{FK}}=\max_{0\leq i,j\leq k}|x_i-x_j|$. This identity shows that $\|\cdot\|_{\mathrm{FK}}$ is a norm.

With this notation, for every $v\in h\Z^k$ and every $x\in\R^k$, we have
\begin{equation}\label{eq:POU-beta-v}
\beta_v(x)=\max\left\{0,1-\frac{1}{h}\|x-v\|_{\mathrm{FK}}\right\}.
\end{equation}
Indeed, a direct verification on every Freudenthal--Kuhn simplex shows that the right-hand side of \eqref{eq:POU-beta-v} is affine on that simplex. Moreover, at every point $w\in h\Z^k$, it equals one if $w=v$ and zero otherwise. Thus \eqref{eq:POU-beta-v} follows from the uniqueness in the definition of $\beta_v$.

We define the coloring $c\colon h\Z^k\to\Z/(k+1)\Z$ by
\[c(v)=\frac{1}{h}\sum_{i=1}^k v_i+(k+1)\Z\]
for every $v=(v_1,\ldots,v_k)\in h\Z^k$. This function is well defined because $v_i/h\in\Z$ for every $1\leq i\leq k$.

\begin{lemma}\label{lemma:beta}
The family $(\beta_v)_{v\in h\Z^k}$ satisfies the following properties:
\begin{enumerate}
\item For every $v,w\in h\Z^k$, if $c(v)=c(w)$ and the functions $\beta_v$ and $\beta_w$ have a common positive point, then $v=w$. In particular, for every $x\in\R^k$, the set $\{v\in h\Z^k:\beta_v(x)>0\}$ consists of at most $k+1$ elements.
\item For every $x\in\R^k$, we have $\sum_{v\in h\Z^k}\beta_v(x)=1$.
\item For every $v\in h\Z^k$, we have
\[\supp\beta_v\subseteq v+[-h,h]^k,\qquad \diam(\supp\beta_v)\leq2h\sqrt{k}.\]
Here the diameter is taken with respect to the Euclidean norm.
\end{enumerate}
\end{lemma}

\begin{proof}
\begin{enumerate}
\item Let $v,w\in h\Z^k$ satisfy $c(v)=c(w)$, and suppose that there exists $x\in\R^k$ such that $\beta_v(x)>0$ and $\beta_w(x)>0$. By \eqref{eq:POU-beta-v}, we have
\[\|x-v\|_{\mathrm{FK}}<h,\qquad \|x-w\|_{\mathrm{FK}}<h.\]
Therefore, the triangle inequality gives $\|v-w\|_{\mathrm{FK}}<2h$. Put $q=(v-w)/h\in\Z^k$. By the homogeneity of the Freudenthal--Kuhn norm, we have $\|q\|_{\mathrm{FK}}<2$. Since the coordinates of $q$ are integers, \eqref{eq:FK-norm} shows that $\|q\|_{\mathrm{FK}}$ is a nonnegative integer. Hence $\|q\|_{\mathrm{FK}}\leq1$.

Write $q=(q_1,\ldots,q_k)$. It follows from \eqref{eq:FK-norm} that either $q_i\in\{0,1\}$ for every $1\leq i\leq k$, or $q_i\in\{0,-1\}$ for every $1\leq i\leq k$. In particular, $\left|\sum_{i=1}^kq_i\right|\leq k$.

On the other hand, the equality $c(v)=c(w)$ implies $\sum_{i=1}^kq_i\in(k+1)\Z$. The only multiple of $k+1$ whose absolute value is at most $k$ is zero. Hence $\sum_{i=1}^kq_i=0$. Since all the coordinates of $q$ are nonnegative or all of them are nonpositive, it follows that $q_i=0$ for every $1\leq i\leq k$. Thus $q=0$, and consequently $v=w$.

For every $x\in\R^k$, the preceding argument shows that there is at most one $v\in h\Z^k$ of each color such that $\beta_v(x)>0$. Since the coloring has $k+1$ possible values, the set $\{v\in h\Z^k:\beta_v(x)>0\}$ consists of at most $k+1$ elements.

\item Let $x\in\R^k$. Since the collection $(\Delta_{m,\sigma})_{m\in\Z^k,\,\sigma\in S_k}$ covers $\R^k$, there exists a simplex $\Delta$ that contains $x$. Denote its vertices by $p_0,\ldots,p_k\in h\Z^k$.

Let $v\in h\Z^k\setminus\{p_0,\ldots,p_k\}$. The function $\beta_v$ vanishes at every vertex of $\Delta$. Since $\beta_v$ is affine on $\Delta$, it vanishes identically on $\Delta$. In particular, $\beta_v(x)=0$. It follows that $\sum_{v\in h\Z^k}\beta_v(x)=\sum_{j=0}^k\beta_{p_j}(x)$.

The function $\sum_{j=0}^k\beta_{p_j}$ is affine on $\Delta$. Moreover, for every $0\leq i\leq k$, we have $\sum_{j=0}^k\beta_{p_j}(p_i)=1$. Since an affine function on a simplex is determined by its values at the vertices, the function $\sum_{j=0}^k\beta_{p_j}$ is identically equal to $1$ on $\Delta$. In particular, $\sum_{v\in h\Z^k}\beta_v(x)=1$, as required.

\item Let $v\in h\Z^k$. By \eqref{eq:POU-beta-v}, for every $x\in\R^k$, we have $\beta_v(x)>0$ if and only if $\|x-v\|_{\mathrm{FK}}<h$. Since the Freudenthal--Kuhn norm is continuous, we have
\[\supp\beta_v\subseteq\left\{x\in\R^k:\|x-v\|_{\mathrm{FK}}\leq h\right\}.\]
Let $x\in\supp\beta_v$. Write $x=(x_1,\ldots,x_k)$ and $v=(v_1,\ldots,v_k)$. By \eqref{eq:FK-norm}, for every $1\leq i\leq k$, we have $|x_i-v_i|\leq\|x-v\|_{\mathrm{FK}}\leq h$. Therefore, $x\in v+[-h,h]^k$. This shows that $\supp\beta_v\subseteq v+[-h,h]^k$.

The Euclidean diameter of the cube $v+[-h,h]^k$ is $2h\sqrt{k}$. Since $\supp\beta_v$ is contained in this cube, we conclude that $\diam(\supp\beta_v)\leq2h\sqrt{k}$. This completes the proof.
\end{enumerate}
\end{proof}

Property \textup{(3)} also implies that the family $(\beta_v)_{v\in h\Z^k}$ is locally finite. Indeed, let $Q\subseteq\R^k$ be compact. If $Q\cap\supp\beta_v\neq\varnothing$, then $v\in Q+[-h,h]^k$. Since this set is compact, it contains only finitely many points of $h\Z^k$.

For every color $a\in\Z/(k+1)\Z$, define the function $H_a^0\colon\R^k\to[0,1]$ by
\[H_a^0=\sum_{\substack{v\in h\Z^k\\c(v)=a}}\beta_v.\]
The local finiteness of $(\beta_v)_{v\in h\Z^k}$ implies that this function is well defined and continuous. Moreover, Lemma~\ref{lemma:beta} implies that, at every point of $\R^k$, at most one summand in this sum is nonzero. For every subset $A\subseteq\Z/(k+1)\Z$, define the function $H_A^0\colon\R^k\to[0,1]$ by $H_A^0=\sum_{a\in A}H_a^0$. These functions are continuous and piecewise affine. Lemma~\ref{lemma:beta} also implies that $(H_0^0,\ldots,H_k^0)$ is a nonnegative partition of unity on $\R^k$.

A cyclic interval in $\Z/(k+1)\Z$ is the empty set, the full group, or a set $\{a,a+1,\ldots,a+\ell-1\}$, where the residues are taken modulo $k+1$ and $1\leq\ell\leq k$.

\begin{claim}
For every color $a\in\Z/(k+1)\Z$, we have
\[\{H_a^0>0\}=\bigsqcup_{\substack{v\in h\Z^k\\c(v)=a}}\{\beta_v>0\}.\]
Moreover, the sets $\{\beta_v>0\}$ appearing in this union are precisely the connected components of $\{H_a^0>0\}$, and
\[\overline{\{\beta_v>0\}}=\supp\beta_v\subseteq v+[-h,h]^k.\]
\end{claim}

\begin{proof}
Lemma~\ref{lemma:beta} implies that, at every point of $\R^k$, at most one function $\beta_v$ of color $a$ is positive. Therefore,
\[\{H_a^0>0\}=\bigsqcup_{\substack{v\in h\Z^k\\c(v)=a}}\{\beta_v>0\}.\]
By \eqref{eq:POU-beta-v}, we have $\{\beta_v>0\}=\{x\in\R^k:\|x-v\|_{\mathrm{FK}}<h\}$. This set is convex and therefore connected. Moreover, its complement in $\{H_a^0>0\}$ is the union of the other open sets appearing in the disjoint union. Thus $\{\beta_v>0\}$ is both open and closed in $\{H_a^0>0\}$, and hence these sets are precisely its connected components. Since $\beta_v$ is nonnegative, its nonzero set is $\{\beta_v>0\}$. Therefore, by the definition of support,
\[\overline{\{\beta_v>0\}}=\supp\beta_v\subseteq v+[-h,h]^k.\]
\end{proof}

\begin{claim}
Let $A\subseteq\Z/(k+1)\Z$ be a cyclic interval. For almost every point $x\in\R^k$, the function $H_A^0$ is differentiable at $x$ and satisfies
\[\|d_xH_A^0\|_1\leq\frac{2}{h}.\]
Moreover, $H_A^0$ is $\frac{2}{h}$-Lipschitz with respect to the $\ell^\infty$ norm.
\end{claim}

\begin{proof}
First, let us show that, for every simplex $\Delta_{m,\sigma}$ from the triangulation, the function $H_A^0$ is differentiable at every point $x\in\Int(\Delta_{m,\sigma})$ and satisfies $\|d_xH_A^0\|_1\leq\frac{2}{h}$. Let $m\in\Z^k$ and $\sigma\in S_k$. The restriction of $H_A^0$ to $\Delta_{m,\sigma}$ is affine. This means that there exist a linear map $L_{m,\sigma,A}\colon\R^k\to\R$ and a constant $b_{m,\sigma,A}\in\R$ such that $H_A^0(x)=L_{m,\sigma,A}(x)+b_{m,\sigma,A}$ for every $x\in\Delta_{m,\sigma}$. In particular, for every $x\in\Int(\Delta_{m,\sigma})$, the function $H_A^0$ is differentiable at $x$ and satisfies $d_xH_A^0=L_{m,\sigma,A}$.

Denote the vertices of $\Delta_{m,\sigma}$ by $p_0,\ldots,p_k$, where $p_0=hm$ and
\[p_j=h\left(m+e_{\sigma(1)}+\cdots+e_{\sigma(j)}\right)\]
for every $1\leq j\leq k$. For every $0\leq j\leq k$, put
\[\epsilon_j=\begin{cases}1,&c(p_j)\in A,\\0,&c(p_j)\notin A,\end{cases}\]
and note that $H_A^0(p_j)=\epsilon_j$. Moreover, for every $1\leq j\leq k$, we have $p_j-p_{j-1}=he_{\sigma(j)}$ and $c(p_j)=c(p_{j-1})+1$ in $\Z/(k+1)\Z$. Thus the colors of $p_0,\ldots,p_k$ are all the elements of $\Z/(k+1)\Z$, listed in their cyclic order. Since $A$ is a cyclic interval, the sequence $\epsilon_0,\ldots,\epsilon_k$ changes its value at most twice. Therefore,
\[\sum_{j=1}^k|\epsilon_j-\epsilon_{j-1}|\leq2.\]
For every $1\leq j\leq k$, we have
\[L_{m,\sigma,A}(e_{\sigma(j)})=\frac{1}{h}L_{m,\sigma,A}(p_j-p_{j-1})=\frac{1}{h}\left(H_A^0(p_j)-H_A^0(p_{j-1})\right)=\frac{1}{h}\left(\epsilon_j-\epsilon_{j-1}\right).\]
Consequently, for every $x\in\Int(\Delta_{m,\sigma})$, we have
\[\|d_xH_A^0\|_1=\|L_{m,\sigma,A}\|_1=\sum_{i=1}^k|L_{m,\sigma,A}(e_i)|=\sum_{j=1}^k|L_{m,\sigma,A}(e_{\sigma(j)})|=\frac{1}{h}\sum_{j=1}^k|\epsilon_j-\epsilon_{j-1}|\leq\frac{2}{h}.\]

The complement of the union of the interiors of the Freudenthal--Kuhn simplices equals the union of their proper faces. Since there are countably many simplices and every proper face has Lebesgue measure zero, this complement has Lebesgue measure zero. Thus $H_A^0$ is differentiable almost everywhere and $\|d_xH_A^0\|_1\leq\frac{2}{h}$ for almost every $x\in\R^k$.

It remains to prove the Lipschitz estimate. Let $x,y\in\R^k$. If $x=y$, we have nothing to prove. Assume that $x\neq y$. Since the triangulation is locally finite, the segment joining $x$ and $y$ meets only finitely many simplices. Therefore, there exist distinct points $x=z_0,z_1,\ldots,z_N=y$ that occur in this order along the segment, such that, for every $1\leq j\leq N$, the segment joining $z_{j-1}$ and $z_j$ is contained in one Freudenthal--Kuhn simplex. Since these points occur in order along the segment, there exist $\lambda_1,\ldots,\lambda_N>0$ such that $z_j-z_{j-1}=\lambda_j(y-x)$ for every $1\leq j\leq N$, and $\sum_{j=1}^N\lambda_j=1$. Consequently,
\[\sum_{j=1}^N\|z_j-z_{j-1}\|_\infty=\sum_{j=1}^N\|\lambda_j(y-x)\|_\infty=\left(\sum_{j=1}^N\lambda_j\right)\|y-x\|_\infty=\|y-x\|_\infty.\]
For every $1\leq j\leq N$, choose a Freudenthal--Kuhn simplex containing the segment joining $z_{j-1}$ and $z_j$. The restriction of $H_A^0$ to this simplex is affine, and the $\ell^1$ norm of its differential is at most $\frac{2}{h}$. Therefore,
\[|H_A^0(z_j)-H_A^0(z_{j-1})|\leq\frac{2}{h}\|z_j-z_{j-1}\|_\infty.\]
Adding these inequalities gives
\[|H_A^0(x)-H_A^0(y)|\leq\frac{2}{h}\sum_{j=1}^N\|z_j-z_{j-1}\|_\infty=\frac{2}{h}\|x-y\|_\infty.\]
This shows that $H_A^0$ is $\frac{2}{h}$-Lipschitz with respect to the $\ell^\infty$ norm.
\end{proof}

Fix a nonnegative compactly supported smooth function $\eta\colon\R^k\to[0,\infty)$ satisfying $\int_{\R^k}\eta(z)\,dz=1$. For every $\tau>0$, define the function $\eta_\tau\colon\R^k\to[0,\infty)$ by $\eta_\tau(y)=\tau^{-k}\eta(y/\tau)$ for every $y\in\R^k$. For every $\tau>0$ and every color $a\in\Z/(k+1)\Z$, define the function $H_a^\tau\colon\R^k\to\R$ by convolution:
\[H_a^\tau(y)=(H_a^0*\eta_\tau)(y)=\int_{\R^k}\eta_\tau(z)H_a^0(y-z)\,dz\]
for every $y\in\R^k$. Finally, for every $\tau>0$ and every subset $A\subseteq\Z/(k+1)\Z$, define the function $H_A^\tau\colon\R^k\to\R$ by $H_A^\tau=\sum_{a\in A}H_a^\tau$.

\begin{prop}\label{prop:smooth-colored-partition}
The functions defined above satisfy the following properties:
\begin{enumerate}
\item For every $\tau>0$, the family $(H_0^\tau,\ldots,H_k^\tau)$ is a smooth nonnegative partition of unity on $\R^k$.
\item For every $\tau>0$, every cyclic interval $A\subseteq\Z/(k+1)\Z$, and every $x\in\R^k$, we have
\[\|d_xH_A^\tau\|_1\leq\frac{2}{h}.\]
\item For every color $a\in\Z/(k+1)\Z$, we have $\lim_{\tau\to0^+}\|H_a^\tau-H_a^0\|_\infty=0$.
\end{enumerate}
\end{prop}

\begin{proof}
\begin{enumerate}
\item Let $\tau>0$. Since $\eta_\tau$ and the functions $H_a^0$ are nonnegative, the functions $H_a^\tau$ are nonnegative. They are smooth by the standard properties of convolution. Moreover,
\[\sum_{a=0}^kH_a^\tau=\sum_{a=0}^k(H_a^0*\eta_\tau)=\left(\sum_{a=0}^kH_a^0\right)*\eta_\tau=1*\eta_\tau=1.\]
This proves the first assertion.

\item Let $\tau>0$, and let $A\subseteq\Z/(k+1)\Z$ be a cyclic interval. Since convolution is linear, we have $H_A^\tau=H_A^0*\eta_\tau$. Let $x,y\in\R^k$. Since $H_A^0$ is $\frac{2}{h}$-Lipschitz with respect to the $\ell^\infty$ norm, we have
\[|H_A^\tau(x)-H_A^\tau(y)|\leq\int_{\R^k}\eta_\tau(z)|H_A^0(x-z)-H_A^0(y-z)|\,dz\leq\frac{2}{h}\|x-y\|_\infty.\]
Thus $H_A^\tau$ is also $\frac{2}{h}$-Lipschitz with respect to the $\ell^\infty$ norm. Since $H_A^\tau$ is smooth, it follows that $\|d_xH_A^\tau\|_1\leq\frac{2}{h}$ for every $x\in\R^k$. This proves the second assertion.

\item Let $\tau>0$ and let $a\in\Z/(k+1)\Z$. Since the singleton $\{a\}$ is a cyclic interval, the function $H_a^0$ is $\frac{2}{h}$-Lipschitz with respect to the $\ell^\infty$ norm. Therefore, for every $x\in\R^k$,
\[|H_a^\tau(x)-H_a^0(x)|\leq\int_{\R^k}\eta_\tau(z)|H_a^0(x-z)-H_a^0(x)|\,dz\leq\frac{2}{h}\int_{\R^k}\eta_\tau(z)\|z\|_\infty\,dz.\]
The change of variables $z=\tau w$ gives $\int_{\R^k}\eta_\tau(z)\|z\|_\infty\,dz=\tau\int_{\R^k}\eta(w)\|w\|_\infty\,dw$. Consequently,
\[\|H_a^\tau-H_a^0\|_\infty\leq\frac{2\tau}{h}\int_{\R^k}\eta(w)\|w\|_\infty\,dw.\]
The integral on the right-hand side is finite because $\eta$ is compactly supported. Letting $\tau\to0^+$ proves the third assertion.
\end{enumerate}
\end{proof}

\subsubsection{The functional and the approximation lemmas}
\label{sec:functional}

In this subsection, we establish some basic properties of a semi-symplectic quasi-state that will be useful in the proof of Theorem~\ref{thm:thick-fiber}. Throughout the remainder of this section, let $(M,\omega)$ be a closed symplectic manifold, let $C>0$, and let $\zeta\colon C(M)\to\R$ be a semi-symplectic quasi-state with constant $C$. We begin with two basic properties that are well known for partial symplectic quasi-states. The proofs in our setting are almost identical.

\begin{lemma}\label{lemma:ssqs_basis_properties}
The semi-symplectic quasi-state $\zeta$ satisfies the following properties:
\begin{enumerate}
\item For every $F\in C(M)$ and every $c\in\R$, we have $\zeta(F+c)=\zeta(F)+c$.
\item For every $F,G\in C(M)$, we have
\[\left|\zeta(F)-\zeta(G)\right|\leq\|F-G\|_\infty.\]
\end{enumerate}
\end{lemma}

\begin{proof}
\begin{enumerate}
\item We first consider the case where $F\in C^\infty(M)$. Let $c\in\R$. Since $\{F,c\}=0$, applying the one-sided Poisson-bracket inequality first to $F$ and $c$, and then to $F+c$ and $-c$, gives
\[\zeta(F+c)\leq\zeta(F)+\zeta(c)=\zeta(F)+c\]
and
\[\zeta(F)\leq\zeta(F+c)+\zeta(-c)=\zeta(F+c)-c.\]
Therefore, $\zeta(F+c)=\zeta(F)+c$ for every $F\in C^\infty(M)$ and every $c\in\R$.

Let us now consider the general case. Let $F\in C(M)$ and $c\in\R$. Let $\epsilon>0$, and choose a smooth function $Q\in C^\infty(M)$ such that $\|Q-F\|_\infty<\epsilon$. We have
\[Q-\epsilon\leq F\leq Q+\epsilon\]
and
\[Q+c-\epsilon\leq F+c\leq Q+c+\epsilon.\]
By monotonicity and the smooth case proved above, we obtain
\[\zeta(Q)-\epsilon\leq\zeta(F)\leq\zeta(Q)+\epsilon\]
and
\[\zeta(Q)+c-\epsilon\leq\zeta(F+c)\leq\zeta(Q)+c+\epsilon.\]
Consequently,
\[\left|\zeta(F+c)-\zeta(F)-c\right|\leq2\epsilon.\]
Since this inequality holds for every $\epsilon>0$, we conclude that $\zeta(F+c)=\zeta(F)+c$.

\item Let $F,G\in C(M)$, and put $d=\|F-G\|_\infty$. We have $G-d\leq F\leq G+d$. By monotonicity and the first item,
\[\zeta(G)-d\leq\zeta(F)\leq\zeta(G)+d.\]
Therefore,
\[\left|\zeta(F)-\zeta(G)\right|\leq d=\|F-G\|_\infty.\]
This completes the proof.
\end{enumerate}
\end{proof}
The following proposition is a version of subadditivity for continuous functions adapted to our setting. This is not the most general form that can be proved.\footnote{Cardin--Viterbo defined in \cite[Definition~1.1]{Cardin_Viterbo_2008_Commuting_Hamiltonians_and_Hamilton-Jacobi_multi-time_equations} continuous Poisson-commuting functions using smooth approximations whose Poisson brackets converge uniformly to zero. The Poisson-bracket inequality and Lemma~\ref{lemma:ssqs_basis_properties} then imply subadditivity. More generally, the author, together with Y.~Ganor, L.~Polterovich, and F.~Zapolsky, defined in \cite[Definition~1.27]{DGPZ_2024_Symp_top_and_IVMs} continuous involutive maps $\rho\colon M\to Y$. The same argument gives subadditivity for functions of the form $f_1\circ\rho,\ldots,f_N\circ\rho$.} However, the following form is sufficient for our purposes.

\begin{prop}[Subadditivity for disjoint positivity sets]\label{prop: subadditivity for fns w disjoint pos set}
Let $N\geq1$ be an integer, and let $u_1,\ldots,u_N\in C(M)$ be nonnegative functions satisfying $u_i(x)u_j(x)=0$ for every $x\in M$ and every distinct $i,j\in\{1,\ldots,N\}$. Then
\[\zeta\left(\sum_{i=1}^N u_i\right)\leq\sum_{i=1}^N\zeta(u_i).\]
\end{prop}

To prove Proposition~\ref{prop: subadditivity for fns w disjoint pos set}, we need the following auxiliary lemma.

\begin{lemma}\label{lemma:positive-approximation}
For every nonnegative continuous function $u\in C(M)$, there exists a sequence $(u_\ell)_{\ell\geq1}$ of nonnegative smooth functions in $C^\infty(M)$ that converges uniformly to $u$ and satisfies $\supp u_\ell\subseteq u^{-1}((0,\infty))$ for every integer $\ell\geq1$.
\end{lemma}

\begin{proof}
Let $\ell\geq1$, and choose a smooth function $q_\ell\in C^\infty(M)$ satisfying $\|q_\ell-u\|_\infty<\frac{1}{2\ell}$. Choose a smooth function $\chi_\ell\colon M\to[0,1]$ that is equal to one on $u^{-1}([2/\ell,\infty))$ and satisfies $\supp\chi_\ell\subseteq u^{-1}((1/\ell,\infty))$. Define the smooth function $u_\ell\colon M\to\R$ by $u_\ell=\chi_\ell \cdot q_\ell$.

Let us first show that $u_\ell$ is nonnegative. Let $x\in M$. If $\chi_\ell(x)=0$, then $u_\ell(x)=0$. If $\chi_\ell(x)>0$, then $x\in\supp\chi_\ell$, and hence $u(x)>1/\ell$. Therefore, 
\[q_\ell(x)>u(x)-\frac{1}{2\ell}>\frac{1}{2\ell}>0.\]
Since $\chi_\ell(x)\geq0$, it follows that $u_\ell(x)\geq0$. Moreover,
\[
\supp u_\ell \subseteq\supp\chi_\ell \subseteq u^{-1}((1/\ell,\infty)) \subseteq u^{-1}((0,\infty)).
\]

Let us now prove that $\|u_\ell-u\|_\infty\leq\frac{5}{2\ell}$. Let $x\in M$. If $u(x)\geq2/\ell$, then $\chi_\ell(x)=1$, and therefore
\[|u_\ell(x)-u(x)|=|q_\ell(x)-u(x)|\leq\|q_\ell-u\|_\infty<\frac{1}{2\ell}.\]
Otherwise, $u(x)<2/\ell$, and hence
\[|u_\ell(x)-u(x)|\leq\chi_\ell(x)|q_\ell(x)-u(x)|+(1-\chi_\ell(x))u(x)<\frac{1}{2\ell}+\frac{2}{\ell}=\frac{5}{2\ell}.\]
This shows that $\|u_\ell-u\|_\infty\leq\frac{5}{2\ell}$. Thus the sequence $(u_\ell)_{\ell\geq1}$ converges uniformly to $u$, as required.
\end{proof}

\begin{proof}[Proof of Proposition~\ref{prop: subadditivity for fns w disjoint pos set}]
For every $i\in\{1,\ldots,N\}$, apply Lemma~\ref{lemma:positive-approximation} to $u_i$ to obtain a sequence $(u_{i,\ell})_{\ell\geq1}$ of nonnegative smooth functions that converges uniformly to $u_i$ and satisfies $\supp u_{i,\ell}\subseteq u_i^{-1}((0,\infty))$ for every integer $\ell\geq1$.

Fix an integer $\ell\geq1$. Since the positivity sets of $u_1,\ldots,u_N$ are pairwise disjoint, the supports of $u_{1,\ell},\ldots,u_{N,\ell}$ are pairwise disjoint. Therefore, $\{u_{i,\ell},u_{j,\ell}\}=0$ for every distinct $i,j\in\{1,\ldots,N\}$. By the bilinearity of the Poisson bracket, the functions $\sum_{i=1}^j u_{i,\ell}$ and $u_{j+1,\ell}$ Poisson commute for every $1\leq j\leq N-1$. Repeated application of the one-sided Poisson-bracket inequality gives
\[\zeta\left(\sum_{i=1}^N u_{i,\ell}\right)\leq\sum_{i=1}^N\zeta(u_{i,\ell}).\]
Since $N$ is finite, $\sum_{i=1}^N u_{i,\ell}$ converges uniformly to $\sum_{i=1}^N u_i$ as $\ell\to\infty$. Therefore, Lemma~\ref{lemma:ssqs_basis_properties} allows us to pass to the limit in the preceding inequality and obtain
\[\zeta\left(\sum_{i=1}^N u_i\right)\leq\sum_{i=1}^N\zeta(u_i).\]
This completes the proof.
\end{proof}

\subsubsection{The pseudoheavy thick-fiber theorem}
\label{sec:thick}

We continue to use the notation introduced above and at the beginning of Section~\ref{ss: proof of thick-fiber}.

\begin{lemma}\label{lemma:chain}
Let $k\geq1$ be an integer, let $F\colon M\to\R^k$ be a smooth map, and let $u,v\colon\R^k\to\R$ be smooth functions satisfying $\|du\|_{1,\infty}<\infty$ and $\|dv\|_{1,\infty}<\infty$. Then
\[\|\{u\circ F,v\circ F\}\|_\infty\leq pb(F)\|du\|_{1,\infty}\|dv\|_{1,\infty}.\]
\end{lemma}

\begin{proof}
Let $x\in M$. By the chain rule and the bilinearity of the Poisson bracket, we have
\[\{u\circ F,v\circ F\}(x)=\sum_{i,j=1}^k\partial_i u(F(x))\partial_j v(F(x))\{F_i,F_j\}(x).\]
Therefore,
\[|\{u\circ F,v\circ F\}(x)|\leq pb(F)\left(\sum_{i=1}^k|\partial_i u(F(x))|\right)\left(\sum_{j=1}^k|\partial_j v(F(x))|\right)\leq pb(F)\|du\|_{1,\infty}\|dv\|_{1,\infty}.\]
Thus, we deduce that 
\[\|\{u\circ F,v\circ F\}\|_\infty\leq pb(F)\|du\|_{1,\infty}\|dv\|_{1,\infty},\]
as required.
\end{proof}

For every $h>0$ and $k\in\N$, consider the partition of unity $(\beta_v)_{v\in h\Z^k}$ constructed in Section~\ref{ss:partition}. We also use the color sums $H_a^\tau$, where $a\in\Z/(k+1)\Z$ and $\tau\geq0$, defined there. The argument below is also related to the successive-addition argument used in \cite[Theorem~1.7]{Entov_Polterovich_Zapolsky_2007_quasi-morphisms_and_PB}.

\begin{prop}\label{prop:witness} 
Let $k\in\N$, let $F\colon M\to\R^k$ be a smooth map, and let $h>0$. For every $v\in h\Z^k$, define the continuous function $u_v\colon M\to[0,1]$ by $u_v=\beta_v\circ F$. Then only finitely many functions in the family $(u_v)_{v\in h\Z^k}$ are nonzero, and we have
\begin{equation}\label{eq:witness}
1\leq\sum_{v\in h\Z^k}\zeta(u_v)+\frac{2C(k-1)\sqrt{pb(F)}}{h}.
\end{equation}
\end{prop}\begin{proof}
Put $K=F(M)$ and $\delta=pb(F)$. Since $K$ is compact and the family $(\beta_v)_{v\in h\Z^k}$ is locally finite, only finitely many functions $\beta_v$ have supports that intersect $K$. Therefore, only finitely many functions $(u_v)_{v\in h\Z^k}$ are nonzero.
For every $a\in\{0,\ldots,k\}$, define the continuous function $g_a\colon M\to[0,1]$ by
\[g_a=H_a^0\circ F=\sum_{\substack{v\in h\Z^k\\c(v)=a}}u_v.\]
By Lemma~\ref{lemma:beta}, the positivity sets of the nonzero functions $u_v$ of the same color are pairwise disjoint. Therefore, Proposition~\ref{prop: subadditivity for fns w disjoint pos set} gives
\begin{equation}\label{eq:color-subadditivity}
\zeta(g_a)\leq\sum_{\substack{v\in h\Z^k\\c(v)=a}}\zeta(u_v)
\end{equation}
for every $a\in\{0,\ldots,k\}$.
For every $\tau>0$ and every $a\in\{0,\ldots,k\}$, define the smooth function $g_a^\tau\colon M\to[0,1]$ by $g_a^\tau=H_a^\tau\circ F$. Let us show that, for every $\tau>0$, we have
\begin{equation}\label{eq:desire-ineq-for-witness}
1\leq\sum_{a=0}^k\zeta(g_a^\tau)+\frac{2C(k-1)\sqrt{\delta}}{h}.
\end{equation}
Fix $\tau>0$. For every $m\in\{0,\ldots,k-1\}$, define the smooth function $P_m^\tau\colon M\to[0,1]$ by
\[P_m^\tau=\sum_{a=0}^m g_a^\tau=H_{\{0,\ldots,m\}}^\tau\circ F.\]
We claim that, for every $m\in\{0,\ldots,k-1\}$, we have
\begin{equation}\label{eq:prefix-telescope}
\zeta(P_m^\tau)\leq\sum_{a=0}^m\zeta(g_a^\tau)+\frac{2Cm\sqrt{\delta}}{h}.
\end{equation}
If $k=1$, then \eqref{eq:prefix-telescope} follows immediately from the equality $P_0^\tau=g_0^\tau$. Suppose that $k\geq2$. We prove \eqref{eq:prefix-telescope} by induction on $m\in\{0,\ldots,k-1\}$. For $m=0$, this is an equality because $P_0^\tau=g_0^\tau$. Let $m\in\{1,\ldots,k-1\}$, and assume that
\[\zeta(P_{m-1}^\tau)\leq\sum_{a=0}^{m-1}\zeta(g_a^\tau)+\frac{2C(m-1)\sqrt{\delta}}{h}.\]
We have $P_m^\tau=P_{m-1}^\tau+g_m^\tau$. The sets $\{0,\ldots,m-1\}$ and $\{m\}$ are cyclic intervals in $\Z/(k+1)\Z$. Therefore, Proposition~\ref{prop:smooth-colored-partition} and Lemma~\ref{lemma:chain} imply that
\[\|\{P_{m-1}^\tau,g_m^\tau\}\|_\infty\leq\frac{4\delta}{h^2}.\]
The one-sided Poisson-bracket inequality gives
\[\zeta(P_m^\tau)\leq\zeta(P_{m-1}^\tau)+\zeta(g_m^\tau)+\frac{2C\sqrt{\delta}}{h}.\]
Combining this inequality with the induction assumption, we obtain
\[\zeta(P_m^\tau)\leq\sum_{a=0}^m\zeta(g_a^\tau)+\frac{2Cm\sqrt{\delta}}{h}.\]
This proves \eqref{eq:prefix-telescope}.
Since $(H_0^\tau,\ldots,H_k^\tau)$ is a partition of unity, we have $P_{k-1}^\tau+g_k^\tau=1$. Consequently,
\[\{P_{k-1}^\tau,g_k^\tau\}=\{1-g_k^\tau,g_k^\tau\}=0.\]
By normalization and the one-sided Poisson-bracket inequality,
\[1=\zeta(1)=\zeta(P_{k-1}^\tau+g_k^\tau)\leq\zeta(P_{k-1}^\tau)+\zeta(g_k^\tau).\]
Applying \eqref{eq:prefix-telescope} with $m=k-1$, we obtain
\[1\leq\sum_{a=0}^k\zeta(g_a^\tau)+\frac{2C(k-1)\sqrt{\delta}}{h}.\]
This proves \eqref{eq:desire-ineq-for-witness}.
By Proposition~\ref{prop:smooth-colored-partition}, $H_a^\tau$ converges uniformly to $H_a^0$ as $\tau\to0^+$ for every $a\in\{0,\ldots,k\}$. Therefore, $g_a^\tau$ converges uniformly to $g_a$ as $\tau\to0^+$. Lemma~\ref{lemma:ssqs_basis_properties} allows us to pass to the limit and obtain
\[1\leq\sum_{a=0}^k\zeta(g_a)+\frac{2C(k-1)\sqrt{\delta}}{h}.\]
Summing \eqref{eq:color-subadditivity} over all the colors gives
\[1\leq\sum_{v\in h\Z^k}\zeta(u_v)+\frac{2C(k-1)\sqrt{\delta}}{h}.\]
This proves \eqref{eq:witness}.
\end{proof}

\begin{proof}[Proof of Theorem~\ref{thm:thick-fiber}]
Fix an integer $k\geq1$ and a smooth map $F\colon M\to\R^k$. Put $\delta=pb(F)$ and $R_0=4Ck^{3/2}\sqrt{\delta}$. We prove the required assertion in two steps.
\\

\noindent\textbf{Step 1:} Let us prove that, for every $R>R_0$, there exist a continuous function $w\colon M\to\R$ and a point $p\in F(M)$ such that
\begin{itemize}
\item $\zeta(w)>0$;
\item $\supp w\subseteq F^{-1}(B_p(R))$.
\end{itemize}
Indeed, let $R>R_0=4Ck^{3/2}\sqrt{\delta}$. Since $k-1\leq k$, we have
\[2C(k-1)\sqrt{\delta}\leq2Ck\sqrt{\delta}=\frac{R_0}{2\sqrt{k}}<\frac{R}{2\sqrt{k}}.\]
Therefore, there exists $h>0$ such that
\[2C(k-1)\sqrt{\delta}<h<\frac{R}{2\sqrt{k}}.\]
Consider the partition of unity $(\beta_v)_{v\in h\Z^k}$ of $\R^k$ defined in Subsection~\ref{ss:partition}. For every $v\in h\Z^k$, define the continuous function $u_v\colon M\to\R$ by $u_v=\beta_v\circ F$. By Proposition~\ref{prop:witness}, we have
\[0<1-\frac{2C(k-1)\sqrt{\delta}}{h}\leq\sum_{v\in h\Z^k}\zeta(u_v).\]
This shows that there exists $v\in h\Z^k$ such that $\zeta(u_v)>0$. By monotonicity, there exists $x_0\in M$ such that $u_v(x_0)>0$. Indeed, otherwise $u_v\leq0$, and hence $\zeta(u_v)\leq\zeta(0)=0$, which is a contradiction. Consider $w=u_v$ and $p=F(x_0)$.

From the choice of $w$, we know that $\zeta(w)>0$ and $p\in F(M)$. Moreover, $\beta_v(p)=u_v(x_0)>0$, and hence $p\in\supp\beta_v$. By Lemma~\ref{lemma:beta}, we have
\[\diam(\supp\beta_v)\leq2\sqrt{k}h<R.\]
Since $p\in\supp\beta_v$, this shows that $\supp\beta_v\subseteq B_p(R)$.

Let us show that $\supp w\subseteq F^{-1}(B_p(R))$. Indeed, let $x\in M\setminus F^{-1}(B_p(R))$. Then
\[F(x)\in\R^k\setminus B_p(R)\subseteq\R^k\setminus\supp\beta_v.\]
The set $\R^k\setminus\supp\beta_v$ is open, and $\beta_v$ vanishes on this set. Therefore, $w=\beta_v\circ F$ vanishes on the open neighborhood $F^{-1}(\R^k\setminus\supp\beta_v)$ of $x$. Thus $x\notin\supp w$. This shows that $\supp w\subseteq F^{-1}(B_p(R))$ and completes the proof of Step 1.
\\

\noindent\textbf{Step 2:} We complete the proof of Theorem~\ref{thm:thick-fiber} by proving that there exists $y\in F(M)$ such that $F^{-1}(\overline{B}_y(R_0))$ is $\zeta$-pseudoheavy, and in particular, non-displaceable.

For every integer $j\geq1$, put $R_j=R_0+\frac{1}{j}$. By Step 1, there exist a point $p_j\in F(M)$ and a continuous function $w_j\colon M\to\R$ such that $\zeta(w_j)>0$ and $\supp w_j\subseteq F^{-1}(B_{p_j}(R_j))$.

Since $F(M)$ is compact, the sequence $(p_j)_{j\geq1}$ has a convergent subsequence; denote its limit by $y\in F(M)$. Consider $y$, and let us show that
\[X=F^{-1}\left(\overline{B}_y(R_0)\right)\]
is $\zeta$-pseudoheavy. Indeed, let $U\subseteq M$ be an open neighborhood of $X$, and let us show that there exists a continuous function $w\colon M\to\R$ such that $\supp w\subseteq U$ and $\zeta(w)>0$. If $U=M$, then the constant function $1$ satisfies the required properties. Otherwise, assume that $M\setminus U\neq\varnothing$. Since $M\setminus U$ is compact and is disjoint from $X=F^{-1}\left(\overline{B}_y(R_0)\right)$, we deduce that
\[d=\min_{x\in M\setminus U}\|F(x)-y\|\]
is strictly bigger than $R_0$, and hence
\[F^{-1}\left(B_y(d)\right)\subseteq U.\]
Also, since $(R_j)_{j\geq1}$ converges to $R_0$ and $(p_j)_{j\geq1}$ has a subsequence that converges to $y$, we deduce that there exists an integer $j\in \N$ such that
\[R_j+\|p_j-y\|<d.\]
Consider $w_j$. We know that $\zeta(w_j)>0$ and $\supp w_j\subseteq F^{-1}(B_{p_j}(R_j))$. Note that for every $z\in B_{p_j}(R_j)$, we have
\[\|z-y\|\leq\|z-p_j\|+\|p_j-y\|<R_j+\|p_j-y\|<d,\]
thus $B_{p_j}(R_j)\subseteq B_y(d)$, which implies that
\[\supp w_j\subseteq F^{-1}(B_{p_j}(R_j))\subseteq F^{-1}(B_y(d))\subseteq U.\]
Hence $X=F^{-1}(\overline{B}_y(R_0))$ is $\zeta$-pseudoheavy. Finally, Remark~\ref{rem: pseudoheavy implies non-disp} implies that $F^{-1}(\overline{B}_y(R_0))$ is non-displaceable because it is $\zeta$-pseudoheavy. This completes the proof.
\end{proof}

\bibliographystyle{alpha}
\bibliography{quantitative_thick_fiber}
\end{document}